\documentclass[11pt]{article}
\usepackage[english]{babel}
\usepackage{amsmath,amsthm,amssymb}
\usepackage{mathrsfs}
\usepackage{bbm}
\usepackage{amsfonts}
\usepackage[margin=0.8in]{geometry}
\usepackage{graphicx}
\usepackage{caption}
\usepackage{subcaption}
\usepackage{wasysym}
\usepackage{enumerate}
\usepackage{algorithm2e}
\usepackage{tabularx}
\usepackage{color}

\newtheorem{thm}{Theorem}[section]
\newtheorem{ass}[thm]{Assumption}
\newtheorem{cor}[thm]{Corollary}
\newtheorem{lem}[thm]{Lemma}
\newtheorem{prop}[thm]{Proposition}
\theoremstyle{definition}
\newtheorem{defn}[thm]{Definition}

\theoremstyle{rem}
\newtheorem{rem}[thm]{Remark}
\numberwithin{equation}{section}
\newcommand{\R}{\mathbb R}

\newcommand{\bbF}{\mathbb F}

\newcommand{\mcB}{\mathcal{B}}
\newcommand{\mcD}{\mathcal D}

\newcommand{\mcN}{\mathcal N}
\newcommand{\mcF}{\mathcal F}
\newcommand{\mcI}{\mathcal I}

\newcommand{\mcT}{\mathcal T}
\newcommand{\mcP}{\mathcal P}

\newcommand{\mcU}{\mathcal U}
\newcommand{\mcS}{\mathcal S}
\newcommand{\Sqlc}{\mathcal{S}_{\textit{qlc}}}

\newcommand{\E}{\mathbb{E}}
\newcommand{\Prob}{\mathbb{P}}

\newcommand{\Zed}{\mathbb{Z}}

\newcommand{\vect}{\bold{t}}
\newcommand{\vecb}{\bold{b}}

\newcommand{\esssup}{\mathop{\rm{ess}\,\sup}}

\newcommand{\ett}{\mathbbm{1}}
\newcommand{\bigSET}{\mcD}

\newcommand{\cadlag}{c\`adl\`ag~}

\newcommand{\cadlagp}{c\`adl\`ag.~}
\newcommand{\dimX}{d}

\newcommand{\I}{\mathcal I}
\newcommand{\ie}{\textit{i.e.\ }}
\newcommand{\eg}{\textit{e.g.\ }}

\begin{document}

\title{A Finite Horizon Optimal Switching Problem with Memory and Application to Controlled SDDEs\footnote{This work was supported by the Swedish Energy Authorities through grant number 42982-1}}%: An Application to Hydropower Planning

\author{Magnus Perninge\footnote{M.\ Perninge is with the Department of Physics and Electrical Engineering, Linnaeus University, V\"axj\"o,
Sweden. e-mail: magnus.perninge@lnu.se.}} %
\maketitle
% ----------------------------------------------------------------
\begin{abstract}
We consider an optimal switching problem where the terminal reward depends on the entire control trajectory. We show existence of an optimal control by applying a probabilistic technique based on the concept of Snell envelopes. We then apply this result to solve an impulse control problem for stochastic delay differential equations driven by a Brownian motion and an independent compound Poisson process. Furthermore, we show that the studied problem arises naturally when maximizing the revenue from operation of a group of hydro-power plants with hydrological coupling.
\end{abstract}

% ----------------------------------------------------------------
\section{Introduction}
The standard optimal switching problem (sometimes referred to as starting and stopping problem) is a stochastic optimal control problem of impulse type that arises when an operator controls a dynamical system by switching between the different members in a set of operation modes $\mcI=\{1,\ldots,m\}$. In the two-modes setting ($m=2$) the modes may represent, for example, ``operating'' and ``closed'' when maximizing the revenue from mineral extraction in a mine as in~\cite{Brennan}. In the multi-modes setting the operating modes may represent different levels of power production in a power plant when the owner seeks to maximize her total revenue from producing electricity~\cite{CarmLud} or the states ``operating'' and ``closed'' of single units in a multi-unit production facility as in~\cite{BrekkeOksendal}.

In optimal switching the control takes the form $u=(\tau_1,\ldots,\tau_N;\beta_1,\ldots,\beta_N)$, where $\tau_1\leq\tau_2\leq\cdots\leq\tau_N$ is a sequence of times when the operator intervenes on the system and $\beta_j\in\mcI^{-\beta_{j-1}}:= \mcI\setminus\{\beta_{j-1}\}$ is the mode in which the system is operated during $[\tau_j,\tau_{j+1})$. The standard multi-modes optimal switching problem in finite horizon ($T<\infty$) can be formulated as finding the control that maximizes
\begin{equation*}
\E\bigg[\int_0^T \phi_{\xi_s}(s)ds+\psi_{\xi_T}-\sum_{j=1}^Nc_{\beta_{j-1},\beta_{j}}(\tau_j)\bigg],
\end{equation*}
where $\xi_t=b_0\ett_{[0,\tau_{1})}(t)+\sum_{j=1}^N \beta_j\ett_{[\tau_{j},\tau_{j+1})}(t)$ is the \emph{operation mode} (when starting in a predefined mode $b_0\in\mcI$), $\phi_b$ and $\psi_b$ are the running and terminal reward in mode $b\in\mcI$, respectively and $c_{b,b'}(t)$ is the cost incurred by switching from mode $b$ to mode $b'$ at time $t\in [0,T]$.

The standard optimal switching problem has been thoroughly investigated in the last decades after being popularised in~\cite{Brennan}. In~\cite{HamJean} a solution to the two-modes problem was found by rewriting the problem as an existence and uniqueness problem for a doubly reflected backward stochastic differential equation. In~\cite{BollanMSwitch1} existence of an optimal control for the multi-modes optimal switching problem was shown by a probabilistic method based on the concept of Snell envelopes. Furthermore, existence and uniqueness of viscosity solutions to the related Bellman equation was shown for the case when the switching costs are constant and the underlying uncertainty is modeled by a stochastic differential equation (SDE) driven by a Brownian motion. In~\cite{ElAsri} the existence and uniqueness results of viscosity solutions was extended to the case when the switching costs depend on the state variable. Since then, results have been extended to Knightian uncertainty~\cite{HuTang,HamZhang,ChassElieKharr} and non-Brownian filtration and signed switching costs~\cite{MartyrSigned}. For the case when the underlying uncertainty can be modeled by a diffusion process, generalization to the case when the control enters the drift and volatility term was treated in \cite{EliKharrCV}. This was further developed to include state constraints in~\cite{KharrSC}. Another important generalization is to the case when the operator only has partial information about the present state of the diffusion process as treated in~\cite{UUganget}.

In the present work we consider the setting with running and terminal rewards that depend on the entire history of the control. We also show that a special case of the type of switching problems that we consider is that of a controlled stochastic delay differential equation (SDDE), driven by a finite intensity L\'evy process.

To motivate our problem formulation we consider the situation when an operator of two hydro-power plants, located in the same river, wants to maximize her revenue from producing electricity during a fixed operation period. We assume that each plant has its own water reservoir. The power production in a hydropower plant depends on the drop height from the water level of the reservoir to the outlet and thus on the amount of water in the reservoir. As water that passes through the upstream plant will eventually reach the reservoir of the downstream plant we need to consider part of the control history in the upstream plant when optimizing operation of the downstream plant.

In this setting our cost functional can be written
\begin{align}
J(u)&:=\E\bigg[\int_0^T\phi(s,\tau_1,\ldots,\tau_{N_s};\beta_1,\ldots,\beta_{N_s})ds 
+ \psi(\tau_1,\ldots,\tau_{N};\beta_1,\ldots,\beta_{N}) -\sum_{j}c_{\beta_{j-1},\beta_j}(\tau_j)\bigg],\label{ekv:objfun}
\end{align}
where $N_s:=\max\{j:\tau_j\leq s\}$. The contribution of the present work is twofold. First, we show that the problem of maximizing $J$ can be solved under certain assumptions on $\phi$, $\psi$ and the switching costs $c_{\cdot,\cdot}$ by finding an optimal control in terms of a family of interconnected value processes, that we refer to as a \emph{verification family}. We then show that the revenue maximization problem of the hydro-power producer can be formulated as an impulse control problem where the uncertainty is modeled by a controlled SDDE and use our initial result to find an optimal control for this problem.

The remainder of the article is organized as follows. In the next section we state the problem, set the notation used throughout the article and detail the set of assumptions that are made. Then, in Section~\ref{sec:VERthm} a verification theorem is derived. This verification theorem is an extension of the original verification theorem for the multi-modes optimal switching problem developed in~\cite{BollanMSwitch1} and presumes the existence of a verification family. In Section~\ref{sec:exist} we show that, under the assumptions made, there exists a verification family, thus proving existence of an optimal control for the switching problem with cost functional $J$. In Section~\ref{sec:sdde} we more carefully investigate the example of the hydro-power producer and show that the case of a controlled SDDE fits into the problem description investigated in Sections~\ref{sec:VERthm} and \ref{sec:exist}.

\section{Preliminaries}
We consider a finite horizon problem and thus assume that the terminal time $T$ is fixed with $T<\infty$.\\

We let $(\Omega,\mcF,\bbF,\Prob)$ be a probability space, with $\bbF:=(\mcF_t)_{0\leq t\leq T}$ a filtration satisfying the usual conditions in addition to being quasi-left continuous.

\begin{rem}
Recall here the concept of quasi-left continuity: A \cadlag process $(X_t:0\leq t\leq T)$ is quasi-left continuous if for each predictable stopping time $\gamma$ and every announcing sequence of stopping times $\gamma_k\nearrow\gamma$ we have $X_{\gamma -}:=\lim\limits_{k\to\infty}X_{\gamma_k} = X_\gamma$, $\Prob$-a.s. A filtration is quasi-left continuous if $\mcF_{\gamma}=\mcF_{\gamma-}$ for every predictable stopping time $\gamma$.
\end{rem}

\noindent Throughout we will use the following notation:
\begin{itemize}
  \item $\mcP_{\bbF}$ is the $\sigma$-algebra of $\bbF$-progressively measurable subsets of $[0,T]\times \Omega$.
  \item For $p\geq 1$, we let $\mcS^{p}$ be the set of all $\R$-valued, $\mcP_{\bbF}$-measurable, \cadlag processes $(Z_t: 0\leq t\leq T)$ such that, $\Prob$-a.s., $\E\left[\sup_{t\in[0,T]} |Z_t|^p\right]<\infty$ and let $\Sqlc^{p}$ be the subset of processes that are quasi-left continuous.
%  \item We let $\mcH^{2}$ denote the set of all $\R$-valued $\mcP_{\bbF}$-measurable processes $(Z_t: 0\leq t\leq T)$ such that $\E\left[\int_0^T |Z_t|^2 dt\right]<\infty$.
  \item We let $\mcT$ be the set of all $\bbF$-stopping times and for each $\gamma\in\mcT$ we let $\mcT_\gamma$ be the corresponding subsets of stopping times $\tau$ such that $\tau\geq \gamma$, $\Prob$-a.s.
  \item We let $\mcU$ be the set of all $u=(\tau_1,\ldots,\tau_N;\beta_1,\ldots,\beta_N)$, where $(\tau_j)_{j=1}^N$ is a non-decreasing sequence of $\bbF$-stopping times (such that $\lim_{j\to\infty}\tau_j=T$, $\Prob$-a.s.) and $\beta_j\in\mcI^{-\beta_{j-1}}$ is $\mcF_{\tau_j}$-measurable (with $\beta_0:=b_0$, the initial operation mode).
  \item We let $\mcU^f$ denote the subset of $u\in\mcU$ for which $N$ is finite $\Prob$-a.s.~(\ie\\ $\mcU^f:=\{u\in\mcU:\: \Prob\left[\{\omega\in\Omega : N(\omega)>k, \:\forall k>0\}\right]=0\}$) and for all $k\geq 0$ we let $\mcU^k:=\{u\in\mcU:\:N\leq k\}$. For $\gamma\in\mcT$ we let $\mcU_\gamma$ (and $\mcU_\gamma^f$ resp.~$\mcU_\gamma^k$) be the subset of $\mcU$ (and $\mcU^f$ resp.~$\mcU^k$) with $\tau_1\in\mcT_\gamma$.
  \item We define the set $\mcD:=\{(t_1,\ldots;b_1,\ldots):t_1\leq t_2\leq\cdots,\:b_{j+1}\in\mcI^{-b_j}\}$ and let $\mcD^f$ be the corresponding subset of all finite sequences.% and $\mcD^n$ the subset of sequences of length $n$.
  \item For all $n\geq 0$, we let $\bar\mcI^n:=\{(b_1,\ldots,b_n)\in\mcI^n:\, b_{j}\in\mcI^{-b_{j-1}}\}$ and $\bar\mcT^n:=\{(\eta_1,\ldots,\eta_n)\in\mcT^n:\, \eta_1\leq\eta_2\leq\cdots\leq\eta_n\}$.
  \item For $l\geq 0$, we let $\Pi_l:=\{0,T2^{-l},2T2^{-l},\ldots,T\}$ and define the map $\Gamma^l:\cup_{j\geq 1}\bar\mcT^j \to \cup_{j\geq 1}\bar\mcT^j$ as $\Gamma^l(\eta_1,\ldots,\eta_j):=(\inf\{s\in \Pi_l:\,s \geq\eta_1\},\ldots,\inf\{s\in \Pi_l:\,s\geq \eta_j\})$ for all $\eta\in\bar\mcT^j$.
\end{itemize}
To make notation more efficient we introduce the $\mcF_T$-measurable function:
\begin{align*}%\label{ekv:psidef}
\Psi(\tau_1,\ldots,\tau_N;\beta_1,\ldots,\beta_N)&:=\int_0^T\phi(s,\tau_1,\ldots,\tau_{N_s};\beta_1,\ldots,\beta_{N_s})ds+ \psi(\tau_1,\ldots,\tau_{N};\beta_1,\ldots,\beta_{N}).
\end{align*}

\subsection{Problem formulation}
In the above notation, our problem can be characterized by two objects:
\begin{itemize}
  \item A $\mcF_T\otimes \mcB(\mcD)$-measurable map $\Psi:\mcD \to\R$.
  \item A collection, $(c_{b,b'}:\Omega\times [0,T]\to\R)_{(b,b')\in\bar\mcI^2}$, of $\mcP_{\bbF}$-measurable processes.
\end{itemize}
We will make the following preliminary assumptions on these objects:
\begin{ass}\label{ass:onPSIandC}
\begin{enumerate}[(i)]
  \item\label{ass:psiBND} The function $\Psi$ is $\Prob$-a.s.~right-continuous in the intervention times and bounded in the sense that:
  \begin{enumerate}[a)]
    \item $\sup_{u\in\mcU}\E[ |\Psi(\tau_1,\ldots;\beta_1,\ldots)|^2]<\infty$.
    \item For all $(\vect,\vecb)\in\mcD^f$ and any\footnote{Throughout we will use $t_n$ and $b_n$ to denote that last element in the vector $\vect$ and $\vecb$, respectively, whenever $(\vect,\vecb)\in\mcD^f$.} $b\in\mcI^{-b_n}$ we have\\
     $\sup_{u\in\mcU}\E[ \sup_{s\in[t_n,T]}|\Psi(\vect,s,\tau_1\vee s,\ldots;\vecb,b,\beta_1,\ldots)|^2]<\infty$.
  \end{enumerate}
 % \item\label{ass:psiRC} Given $n\geq 0$ and $\vecb\in\bar\mcI^n$ we assume that for every $n$-tuple of stopping times $\eta:=(\eta_1, \eta_2,\ldots,\eta_n)\in\bar\mcT^n$,
%  \begin{equation*}
%  \lim_{k\to\infty}\sup_{u\in\mcU}\E\big[\sup_{s\in[0,T]}|\Psi(\Gamma^k(\eta),\tau_1,\ldots,\tau_N;\vecb,\beta_1,\ldots,\beta_N) -\Psi(\eta,\tau_1,\ldots,\tau_N;\vecb,\beta_1,\ldots,\beta_N)|^2\big]\to 0.
%  \end{equation*}
  \item\label{ass:@end} For each $(\vect,\vecb)\in\mcD^f$ and any $b\in\mcI^{-b_n}$ we have $\Psi(\vect;\vecb)>\Psi(\vect,T;\vecb,b)-c_{b_n,b}(T)$, $\Prob$-a.s.
  \item\label{ass:onc} We assume that $(c_{b,b'})_{(b,b')\in \bar\mcI^2}\in (\Sqlc^2)^{m(m-1)}$ are such that:
  \begin{enumerate}[a)]
    \item\label{ass:oncPOS} $c_{b,b'}\geq 0$, $\Prob$-a.s.
    \item\label{ass:oncNFL} There is an $\epsilon>0$ such that for each $(t_1,\ldots,t_{n},b_1,\ldots,b_n)$ with $0\leq t_1\leq\cdots\leq t_n\leq T$ and $b_1\in\mcI^{-b_n}$, and $b_j\in\mcI^{-b_{j-1}}$ for $j=2,\ldots,n$, we have
    \begin{equation*}
      c_{b_1,b_2}(t_1)+\cdots+c_{b_n,b_1}(t_{n})\geq \epsilon,
    \end{equation*}
    $\Prob$-a.s.
  \end{enumerate}
\end{enumerate}
\end{ass}

The above assumptions are mainly standard assumptions for optimal switching problems translated to our setting. Assumptions (\ref{ass:psiBND}.a) and (\ref{ass:onc}.a) together imply that the expected maximal reward is finite. Assumption~(\ref{ass:@end}) implies that it is never optimal to switch at the terminal time. We show below that the ``no-free-loop'' condition (\ref{ass:onc}.b) together with (\ref{ass:psiBND}.a) implies that, with probability one, the optimal control (whenever it exists) can only make a finite number of switches.\\

%Each control $u=(\tau_1,\ldots,\tau_N;\beta_1,\ldots,\beta_N)$ defines the switching mode starting in $\vecb\in\mcI$ which is a process $(\xi^{\vecb}_t: 0\leq t\leq T)$ given by
%\begin{equation*}
%\xi^{\vecb}_t=\vecb\ett_{[0,\tau_1)}(t)+\sum_{j=1}^N \beta_j\ett_{[\tau_{j},\tau_{j+1})}(t),
%\end{equation*}
%with $\tau_{N+1}=\infty$. For notational simplicity we write $\xi$ for $\xi^{0}$.

We consider the following problem:\\

\noindent\textbf{Problem 1.} Find $u^*\in\mcU$, such that
\begin{equation}\label{ekv:OPTprob}
J(u^*)=\sup_{u\in\mcU} J(u).
\end{equation}
\qed

\bigskip

%\begin{rem}\label{rem:pos}
%Note that we have
%\begin{align*}
%J(u)=&\:\E\bigg[\Psi(u) -\inf_{(\vect,\vecb)\in\mcD}\Psi(\vect,\vecb) -\sum_{j}c_{\beta_{j-1},\beta_j}(\tau_j)\bigg]+\E\bigg[\inf_{(\vect,\vecb)\in\mcD}\Psi(\vect,\vecb) \bigg].
%\end{align*}
%Hence, we can without loss of generality assume that $\Psi$ is non-negative.
%\end{rem}

As a step in solving Problem 1 we need the following proposition which is a standard result for optimal switching problems and is due to the ``no-free-loop'' condition.
\begin{prop}\label{prop:finSTRAT} Suppose that there is a $u^*\in\mcU$ such that $J(u^*)\geq J(u)$ for all $u\in\mcU$. Then $u^*\in\mcU^f$.
\end{prop}

\noindent\emph{Proof.} Pick $\hat u:=(\hat\tau_1,\ldots,\hat\tau_{\hat N};\hat\beta_1,\ldots,\hat\beta_{\hat N})\in \mcU\setminus \mcU^f$ and let $B:=\{\omega \in\Omega: \hat N(\omega)>k, \:\forall k>0\}$, then $\Prob[B]>0$. Furthermore, if $B$ holds then the switching mode $\xi$ must make an infinite number of loops and
\begin{align*}
J(\hat u)&\leq \sup_{u\in\mcU} \E\big[ |\Psi(\tau_1,\ldots;\beta_1,\ldots)|\big]-\frac{k-m}{m}\epsilon \Prob[B]\leq C-\frac{k}{m}\epsilon \Prob[B],
\end{align*}
for all $k\geq 0$, by Assumptions~\ref{ass:prelim}.(\ref{ass:onc}.b) and \ref{ass:prelim}.(\ref{ass:psiBND}.a). However, again by Assumption~\ref{ass:prelim}.(\ref{ass:psiBND}.a) we have\footnote{Throughout $C$ will denote a generic positive constant that may change value from line to line.} $J(\emptyset)\geq -C$. Hence, $\hat u$ is dominated by the strategy of doing nothing and the assertion follows.\qed

%%%%%%%%%%%%%%%%%%%%%%%%%%%%%%%%%%%%%%%%%%%%%%%%%%%%%%%%%%%%%%%%%%%%%%%%%%%%%%%%%%%%%%%%%%%%%%%%%%%%%%%%%%%%%%%%%%%%%%%%%%%%%%%%%%%%%%%%%%%%%%%%%%

\subsection{The Snell envelope}
In this section we gather the main results concerning the Snell envelope that will be useful later on. Recall that a progressively measurable process $U$ is of class [D] if the set of random variables $\{U_\tau:\tau\in\mcT\}$ is uniformly integrable.

\begin{thm}[The Snell envelope]\label{thm:Snell}
Let $U=(U_t)_{0\leq t\leq T}$ be an $\bbF$-adapted, $\R$-valued, \cadlag process of class [D]. Then there exists a unique (up to indistinguishability), $\R$-valued \cadlag process $Z=(Z_t)_{0\leq t\leq T}$ called the Snell envelope, such that $Z$ is the smallest supermartingale that dominates $U$. Moreover, the following holds (with $\Delta U_t:=U_{t}-U_{t-}$):
\begin{enumerate}[(i)]
  \item\label{Snell:sup} For any stopping time $\gamma$,
    \begin{equation}\label{ekv:SnellZ}
      Z_{\gamma}=\esssup_{\tau\in \mcT_{\gamma}}\E\left[U_\tau\big|\mcF_\gamma\right].
    \end{equation}
  \item\label{Snell:DoobMeyer} The Doob-Meyer decomposition of the supermartingale $Z$ implies the existence of a triple $(M,K^c,K^d)$ where $(M_t:0\leq t\leq T)$ is a uniformly integrable right-continuous martingale, $(K^c_t:0\leq t\leq T)$ is a non-decreasing, predictable, continuous process with $K^c_0=0$ and $(K^d_t:0\leq t\leq T)$ is non-decreasing purely discontinuous predictable with $K^d_0=0$, such that
      \begin{equation}\label{ekv:DoobMeyerDec}
        Z_t=M_t-K^c_t-K^d_t.
      \end{equation}
      Furthermore, $\{\Delta_t K^d>0\}\subset \{\Delta_t U<0\}\cap\{Z_{t-}=U_{t-}\}$ for all $t\in[0,T]$.
  \item\label{Snell:att} Let $\theta\in\mcT$ be given and assume that for any predictable $\gamma\in\mcT_\theta$ and any increasing sequence $\{\gamma_k\}_{k\geq 0}$ with $\gamma_k\in\mcT_\theta$ and $\lim_{k\to\infty}\gamma_k=\gamma$, $\Prob$-a.s, we have %$\E[U-_{\gamma_k}]<\infty$ and
$\limsup_{k\to\infty}U_{\gamma_k}\leq U_{\gamma}$, $\Prob$-a.s. Then, the stopping time $\tau^*_{\theta}$ defined by $\tau^*_{\theta}:=\inf\{s\geq\theta:Z_s=U_s\}\wedge T$ is optimal after $\theta$, \ie
    \begin{equation*}
      Z_{\theta}=\E\left[U_{\tau^*_\theta}\big|\mcF_\theta\right].
    \end{equation*}
    Furthermore, in this setting the Snell envelope, $Z$, is quasi-left continuous, \ie $K^d\equiv 0$.
  \item\label{Snell:lim} Let $U^k$ be a sequence of \cadlag processes converging pointwisely to a \cadlag process $U$ and let $Z^k$ be the Snell envelope of $U^k$. Then the sequence $Z^k$ converges pointwisely to a process $Z$ and $Z$ is the Snell envelope of $U$.
\end{enumerate}
\end{thm}
In the above theorem (\ref{Snell:sup})-(\ref{Snell:att}) are standard. Proofs can be found in \cite{ElKarouiLN} (see \cite{Latifa} for an English version), Appendix D in~\cite{KarShreve2}, \cite{HamRefBSDE} and in the appendix of~\cite{CvitKar}. Statement (\ref{Snell:lim}) was proved in \cite{BollanMSwitch1}.\\

The Snell envelope will be the main tool in showing that Problem 1 has a solution.

%%%%%%%%%%%%%%%%%%%%%%%%%%%%%%%%%%%%%%%%%%%%%%%%%%%%%%%%%%%%%%%%%%%%%%%%%%%%%%%%%%%%%%%%%%%%%%%%%%%%%%%%%%%%%%%%%%%%%%%%%%%%%%%%%%%%%%%%%%%%%%%%%%

\subsection{Additional assumptions on regularity}
From the definition of the Snell envelope it is clear that we need to make some further assumptions on the regularity of the involved processes. To facilitate this we define, for each $(\vect,\vecb)=(t_1,\ldots,t_n;b_1,\ldots,b_n)\in\mcD^f$, the value process corresponding to the control $u\in\mcU$ as
\begin{align*}
V^{\vect;\vecb,u}_s&:=\E\big[\Psi(\vect,t_n\vee s\vee\tau_1,\ldots,t_n\vee s\vee\tau_N;\vecb,\beta_1,\ldots,\beta_N)
- \sum_{j=1}^Nc_{\beta_{j-1},\beta_j}(t_n\vee s\vee\tau_j)|\mcF_s\big],
\end{align*}
with $\beta_0:=b_n$.\\

We make the following additional assumptions:
\begin{ass}\label{ass:prelim}
\begin{enumerate}[(i)]
  \item\label{ass:psiRCint} For each $n\geq 0$ and each $(\eta,\vecb)\in\bar\mcT^n\times\bar\mcI^n$ and $b\in\mcI^{-b_n}$ there is a sequence of maps $(\mcU\to\mcU:u\to \hat u^l)_{l\geq 0}$ such that
    \begin{align*}
        \lim_{l\to\infty}\sup_{u\in\mcU}\E\Big[&\sup_{s\in[0,T]}|(V^{\eta;\vecb,u}_s-V^{\Gamma^l(\eta);\vecb,\hat u^l}_s)^++(V^{\eta,s\vee\eta_n;\vecb,b,u}_s-V^{\Gamma^l(\eta),s\vee\Gamma^l(\eta_n);\vecb,b,\hat u^l}_s)^+|^2\Big]=0.
    \end{align*}
    Furthermore, we have
    \begin{align*}
        \lim_{l\to\infty}\sup_{u\in\mcU_{\Gamma^l(\eta_n)}}\E\Big[&\sup_{s\in[0,T]}|(V^{\Gamma^l(\eta);\vecb,u}_s- V^{\eta;\vecb,u}_s)^+(V^{\Gamma^l(\eta),s\vee\Gamma^l(\eta_n);\vecb,b,u}_s- V^{\eta,s\vee\eta_n;\vecb,b,u}_s)^+|^2\Big]=0.%\label{ekv:UkmU}
    \end{align*}
  \item\label{ass:rBNDYbnd} For all $(\vect,\vecb)\in\mcD^f$ and all $b\in \mcI^{-b_n}$, the process $(\esssup_{u\in\mcU^k}V^{\vect,s\vee t_n;\vecb,b,u}_s:0\leq s\leq T)$ is in $\Sqlc^2$ for $k=0,1,\ldots$
\end{enumerate}
\end{ass}

%%%%%%%%%%%%%%%%%%%%%%%%%%%%%%%%%%%%%%%%%%%%%%%%%%%%%%%%%%%%%%%%%%%%%%%%%%%%%%%%%%%%%%%%%%%%%%%%%%%%%%%%%%%%%%%%%%%%%%%%%%%%%%%%%%%%%%%%%%%%%%%%%%
%%%%%%%%%%%%%%%%%%%%%%%%%%%%%%%%%%%%%%%%%%%%%%%%%%%%%%%%%%%%%%%%%%%%%%%%%%%%%%%%%%%%%%%%%%%%%%%%%%%%%%%%%%%%%%%%%%%%%%%%%%%%%%%%%%%%%%%%%%%%%%%%%%
%%%%%%%%%%%%%%%%%%%%%%%%%%%%%%%%%%%%%%%%%%%%%%%%%%%%%%%%%%%%%%%%%%%%%%%%%%%%%%%%%%%%%%%%%%%%%%%%%%%%%%%%%%%%%%%%%%%%%%%%%%%%%%%%%%%%%%%%%%%%%%%%%%

\section{A verification theorem\label{sec:VERthm}}
The method for solving Problem 1 will be based on deriving an optimal control under the assumption that a specific family of processes exists, and then showing that the family indeed does exist. We will refer to any such family of processes as a \emph{verification family}.
\begin{defn}\label{def:vFAM}
We define a \emph{verification family} to be a family of \cadlag supermartingales $((Y^{\vect;\vecb}_s)_{0\leq s\leq T}: (\vect,\vecb)\in \mcD^f)$ such that:
\begin{enumerate}[a)]
  \item\label{vfass:recur} The family satisfies the recursion
  \begin{align}
  Y^{\vect;\vecb}_s&=\esssup_{\tau \in \mcT_{s\vee t_n}} \E\Big[\ett_{[\tau \geq T]}\Psi(\vect;\vecb)+\ett_{[\tau < T]}\max_{\beta\in\mcI^{-b_n}}\left\{-c_{b_n,\beta}(\tau)+Y^{\vect,\tau;\vecb,\beta}_\tau\right\}\Big| \mcF_s\Big].\label{ekv:Ydef}
  \end{align}
  \item\label{vfass:Ybnd} The family is bounded in the sense that $\sup\limits_{u\in\mcU}\E[\sup\limits_{s\in[0,T]}|Y^{\tau_1,\ldots,\tau_N;\beta_1,\ldots,\beta_N}_s|^2]<\infty$.
  \item\label{vfass:rcINt} For all $n\geq 1$ we have that for every $\vecb\in\bar\mcI^n$ and $\eta\in\bar\mcT^n$,
  \begin{equation}\label{ekv:rcINt1}
  \lim_{l\to\infty}\E\big[\sup_{s\in[0,T]}|Y^{\Gamma^l(\eta);\vecb}_s-Y^{\eta;\vecb}_s|^2\big]= 0
  \end{equation}
  and for all $b\in\mcI^{-b_n}$ we have
  \begin{equation}\label{ekv:rcINt2}
  \lim_{l\to\infty}\E\big[\sup_{s\in[0,T]}|Y^{\Gamma^l(\eta),s\vee\Gamma^l(\eta_n);\vecb,b}_s-Y^{\eta,s\vee\eta_n;\vecb,b}_s|^2\big]= 0.
  \end{equation}
  \item\label{vfass:lce} For every $(\vect,\vecb)\in \mcD^f$ and every $b\in\I^{-b_n}$, the process $(Y^{\vect,s;\vecb,b}_s:0\leq s\leq T)$ is in $\Sqlc^2$.
\end{enumerate}
\end{defn}

The purpose of the present section is to reduce the solution of Problem 1 to showing existence of a verification family. This is done in the following verification theorem:
\begin{thm}\label{thm:vfc}
Assume that there exists a verification family $((Y^{\vect;\vecb}_s)_{0\leq s\leq T}: (\vect,\vecb)\in \mcD^f)$. Then the family is unique (\ie there is at most one verification family, up to indistinguishability) and:
\begin{enumerate}[(i)]
  \item Satisfies $Y_0=\sup_{u\in \mcU} J(u)$ (where $Y:=Y^{\emptyset}$).
  \item Defines the optimal control, $u^*=(\tau_1^*,\ldots,\tau_{N^*}^*;\beta_1^*,\ldots,\beta_{N^*}^*)$, for Problem 1, where $(\tau_j^*)_{1\leq j\leq {N^*}}$ is a sequence of $\bbF$-stopping times given by
  \begin{align*}
  \tau^*_j:=\inf \Big\{&s \geq \tau^*_{j-1}:\:Y_s^{\tau^*_{1},\ldots,\tau^*_{j-1};\beta^*_{1},\ldots,\beta^*_{j-1}}=\max_{\beta\in \mcI^{-\beta^*_{j-1}}}\Big\{-c_{\beta_{j-1}^*,\beta}(s)+ Y^{\tau^*_{1},\ldots,\tau^*_{j-1},s;\beta^*_{1},\ldots,\beta^*_{j-1},\beta}_s\Big\}\Big\}\wedge T,%\label{ekv:taujDEF}
  \end{align*}
  $(\beta_j^*)_{1\leq j\leq {N^*}}$ is defined as a measurable selection of
  \begin{equation*}
  \beta^*_j\in\mathop{\arg\max}_{\beta\in \mcI^{-\beta^*_{j-1}}}\Big\{-c_{\beta_{j-1}^*,\beta}(\tau_j^*)+ Y^{\tau^*_{1},\ldots,\tau^*_j;\beta^*_{1},\ldots,\beta^*_{j-1},\beta}_{\tau^*_j}\Big\}
  \end{equation*}
  and $N^*=\max\{j:\tau^*_j<T\}$, with $(\tau_0^*,\beta^*_0):=(0,b_0)$.
\end{enumerate}
\end{thm}

%%%%%%%%%%%%%%%%%%%%%%%%%%%%%%%%%%%%%%%%%%%%%%%%%%%%%%%%%%%%%%%%%%%%%%%%%%%%%%%%%%%%%%%%%%%%%%%%%%%%%%%%%%%%%%%%%%%%%%%%%%%%%%%%%%%%%%%%%%%%%%%%

\noindent\emph{Proof.} The proof is divided into three steps where we first, in steps 1 and 2, show that for any $0\leq j\leq N^*$ we have
\begin{align}\nonumber
&Y^{\tau^*_1,\ldots,\tau^*_j;\beta^*_1,\ldots,\beta^*_j}_{s}=\esssup_{\tau\in\mcT_s}\E\Big[\ett_{[\tau \geq T]}\Psi(\tau^*_1,\ldots,\tau^*_j;\beta^*_1,\ldots,\beta^*_j)
\\ \nonumber
&\qquad+\ett_{[\tau < T]}\max_{\beta\in\mcI^{-\beta^*_j}}\left\{-c_{\beta^*_j,\beta}(\tau) +Y^{\tau^*_1,\ldots,\tau^*_{j},\tau;\beta^*_1,\ldots,\beta^*_{j},\beta}_{\tau}\right\}\Big| \mcF_{s}\Big]
\\ \nonumber
&=\E\Big[\ett_{[\tau^*_{j+1} \geq T]}\Psi(\tau^*_1,\ldots,\tau^*_j;\beta^*_1,\ldots,\beta^*_j)
\\
&\qquad+\ett_{[\tau^*_{j+1} < T]}\left\{-c_{\beta^*_j,\beta^*_{j+1}}(\tau^*_{j+1}) +Y^{\tau^*_1,\ldots,\tau^*_{j+1};\beta^*_1,\ldots,\beta^*_{j+1}}_{\tau^*_{j+1}}\right\}\Big| \mcF_{s}\Big],\label{ekv:esssupATT}
\end{align}
$\Prob$-a.s.~for $s\in[\tau_j^*,\tau^*_{j+1}]$. Then in Step 3 we show that $u^*$ is the optimal control estabilishing \emph{(i)} and \emph{(ii)}. A straightforward generalization to arbitrary initial conditions $(\vect,\vecb)\in\mcD^f$ then gives that
\begin{align}
Y^{\vect;\vecb}_s=\esssup_{u\in \mcU_{s\vee t_n}} \E\Big[\Psi(\vect,\tau_1,\ldots,\tau_N;\vecb,\beta_1,\ldots,\beta_N)-\sum_{j=1}^Nc_{\beta_{j-1},\beta_j}(\tau_j)\Big|\mcF_s\Big],\label{ekv:Yuni}
\end{align}
by which uniqueness follows.\\

\noindent {\bf Step 1} We start by showing that for each $(\vect,\vecb)\in \mcD^f$ the recursion \eqref{ekv:Ydef} can be written in terms of a $\bbF$-stopping time. From \eqref{ekv:Ydef} we note that, by definition, $Y^{\vect;\vecb}$ is the smallest supermartingale that dominates
\begin{align}
U^{\vect;\vecb}:=\Big(&\ett_{[s=T]}\Psi(\vect;\vecb)+\ett_{[s < T]}\max_{\beta\in\mcI^{-b_n}}\big\{-c_{b_n,\beta}(s\vee t_n)+Y^{\vect,s\vee t_n;\vecb,\beta}_s\big\}\Big| :\: 0\leq s\leq T\Big).\label{ekv:dominated}
\end{align}
Now, by Assumption~\ref{ass:onPSIandC}.(\ref{ass:onc}) and property \ref{vfass:lce}) in the definition of a verification family (Definition~\ref{def:vFAM}) we note that $U^{\vect;\vecb}$ is a \cadlag process of class [D] that is quasi-left continuous on $[0,T)$. Furthermore, by Assumption~\ref{ass:onPSIandC}.(\ref{ass:@end}) and property \ref{vfass:lce}) we get that for any sequence $(\eta_k)_{k\geq 0}\subset\mcT$ such that $\eta_k\nearrow T$, $\Prob$-a.s.~we have $\lim_{k\to\infty} U^{\vect;\vecb}_{\eta_k}\leq U^{\vect;\vecb}_T$, $\Prob$-a.s. By Theorem \ref{thm:Snell}.(\ref{Snell:att}) it thus follows that for any $\theta\in \mcT$, there is a stopping time $\gamma_\theta\in \mcT_{t_n\vee\theta}$ such that:
\begin{align*}
Y^{\vect;\vecb}_\theta=\E\Big[\ett_{[\gamma_\theta = T]}\Psi(\vect;\vecb)+\ett_{[\gamma_\theta < T]}\max_{\beta\in\mcI^{-b_n}}\left\{-c_{b_n,\beta}(\gamma_\theta)+Y^{\vect,\gamma_\theta;\vecb,\beta}_{\gamma_\theta}\right\}\Big| \mcF_\theta\Big].
\end{align*}

\noindent {\bf Step 2} We now show that $Y_0=J(u^*)$. We start by noting that $Y$ is the Snell envelope of
\begin{align*}
\Big(\ett_{[s=T]}\Psi_0+\ett_{[s < T]}\max_{\beta\in\mcI^{-b_0}}\left\{-c_{b_0,\beta}(s)+Y^{s,\beta}_s\right\} :\: 0\leq s\leq T\Big),
\end{align*}
where $\Psi_0:=\Psi(\emptyset)$, and by step 1 we thus have
\begin{align*}
Y_0&=\sup_{\tau \in \mcT} \E\Big[\ett_{[\tau=T]}\Psi_0+\ett_{[\tau < T]}\max_{\beta\in\mcI^{-b_0}}\left\{-c_{b_0,\beta}(\tau)+Y^{\tau,\beta}_\tau\right\}\Big]
\\
&=\E\Big[\ett_{[\tau^*_1=T]}\Psi_0+\ett_{[\tau^*_1 < T]}\max_{\beta\in\mcI^{-b_0}}\left\{-c_{b_0,\beta}(\tau^*_1)+Y^{\tau^*_1,\beta}_{\tau^*_1}\right\}\Big]
\\
&=\E\Big[\ett_{[\tau^*_1=T]}\Psi_0+\ett_{[\tau^*_1 < T]}\left\{-c_{b_0,\beta^*_1}(\tau^*_1)+Y^{\tau^*_1,\beta^*_1}_{\tau^*_1}\right\}\Big].
\end{align*}

Moving on we pick $j\in\{1,\ldots, N^*\}$. For $M\geq 0$, let $z_{-1}=-1$ and $z_k:=kT/2^M$ for $k=0,\ldots,2^M$. Furthermore, we define the processes $(\hat Y^{M}_s:0\leq s\leq T)$ and $(\hat U^{M}_t:0\leq s\leq T)$ by
\begin{align*}
\hat Y^{M}_{s}&:=\sum_{(k_1,\ldots k_j)\in \bar\Zed^j}\sum_{(b_1,\ldots,b_j)\in\bar\mcI^j} \E\big[\ett_{(z_{k_1-1},z_{k_1}]}(\tau^*_1)\cdots \ett_{(z_{k_j-1},z_{k_j}]}(\tau^*_j) \ett_{[\beta^*_{1}=b_1]}
\cdots\ett_{[\beta^*_{j}=b_j]}\big|\mcF_s\big] Y_s^{z_{k_1},\ldots,z_{k_{j}};b_1,\ldots,b_j},
\end{align*}
and
\begin{align*}
\hat U^{M}_{s}&:=\sum_{(k_1,\ldots k_j)\in \bar\Zed^j}\sum_{(b_1,\ldots,b_j)\in\bar\mcI^j} \E\big[\ett_{(z_{k_1-1},z_{k_1}]}(\tau^*_1)\cdots \ett_{(z_{k_j-1},z_{k_j}]}(\tau^*_j) \ett_{[\beta^*_{1}=b_1]}\cdots\ett_{[\beta^*_{j}=b_j]}\big|\mcF_s\big]\Big(
\\
&\qquad \ett_{[s=T]}\Psi(z_{k_1},\ldots,z_{k_{j}};b_1,\ldots,b_j)
+\ett_{[s < T]}\max_{\beta\in\mcI^{-b_j}}\Big\{-c_{b_j,\beta}(s\vee z_{k_{j}}) + Y^{z_{k_1},\ldots,z_{k_{j}},s\vee z_{k_{j}};b_1,\ldots,b_j,\beta}_s\Big\}\Big),
\end{align*}
for all $s\in[0,T]$, where $\bar\Zed^j:=\{(k_1,\ldots,k_j)\in\{0,\ldots,2^M\}^j:k_1\leq k_2\leq\cdots\leq k_j\}$. Now, for each $(k_1,\ldots,k_j,b_1,\ldots,b_j)\in \bar\Zed^j\times \bar\mcI^j$ we have that
\begin{align*}
\ett_{(z_{k_1-1},z_{k_1}]}(\tau^*_1)\cdots \ett_{(z_{k_j-1},z_{k_j}]}(\tau^*_j) \ett_{[\beta^*_{1}=b_1]}\cdots \ett_{[\beta^*_{j}=b_j]} Y_s^{z_{k_1},\ldots,z_{k_{j}};b_1,\ldots,b_j},
\end{align*}
is the product of an $\mcF_{\tau^*_j}$--measurable positive r.v.~and a \cadlag supermartingale, thus, it is a \cadlag supermartingale for $s\geq \tau^*_j$. Hence, $\hat Y^{M}$ is the sum of a finite number of \cadlag supermartingales and thus a \cadlag supermartingale itself. By definition we find that $\hat Y^{M}$ dominates $\hat U^{M}$ which is of class [D] by Assumption~\ref{ass:prelim}.(\ref{ass:psiBND}) and property \ref{vfass:Ybnd}). To show that $\hat Y^{M}$ is in fact the Snell envelope of $\hat U^{M}$ assume that $Z$ is another \cadlag supermartingale that dominates $\hat U^{M}$ for all $s\in[\tau^*_{j},T]$. Then for each $(k_1,\ldots,k_j;b_1,\ldots,b_j)\in \bar\Zed^j\times \bar\mcI^j$ and $s\geq \tau_j^*$, we have
\begin{align*}
&\ett_{(z_{k_1-1},z_{k_1}]}(\tau^*_1)\cdots \ett_{(z_{k_j-1},z_{k_j}]}(\tau^*_j) \ett_{[\beta^*_{1}=b_1]}\cdots\ett_{[\beta^*_{j}=b_j]}Z_s
\\
&\geq \ett_{(z_{k_1-1},z_{k_1}]}(\tau^*_1)\cdots \ett_{(z_{k_j-1},z_{k_j}]}(\tau^*_j) \ett_{[\beta^*_{1}=b_1]}\cdots\ett_{[\beta^*_{j}=b_j]}\Big(\Psi(z_{k_1},\ldots,z_{k_j};b_1,\ldots,b_j)
\\
&\qquad+\ett_{[s < T]}\max_{\beta\in\mcI^{-b_j}}\left\{-c_{b_j,\beta}(s)+Y^{z_{k_1},\ldots,z_{k_j},s;b_1,\ldots,b_j,\beta}_s\right\}\Big),
\end{align*}
$\Prob$-a.s.~which by \eqref{ekv:Ydef} gives that
\begin{align*}
&\ett_{(z_{k_1-1},z_{k_1}]}(\tau^*_1)\cdots \ett_{(z_{k_j-1},z_{k_j}]}(\tau^*_j) \ett_{[\beta^*_{1}=b_1]}\cdots\ett_{[\beta^*_{j}=b_j]}Z_s
\\
&\geq \ett_{(z_{k_1-1},z_{k_1}]}(\tau^*_1)\cdots \ett_{(z_{k_j-1},z_{k_j}]}(\tau^*_j) \ett_{[\beta^*_{1}=b_1]}\cdots\ett_{[\beta^*_{j}=b_j]}\hat Y^{z_{k_1},\ldots,z_{k_j};b_1,\ldots,b_j}_s.
\end{align*}
Summing over all $(k_1,\ldots,k_j;b_1,\ldots,b_j)\in \bar\Zed^j\times \bar\mcI^j$ we get $Z_s\geq \hat Y^{M}_s$, $\Prob$-a.s.

Noting that $\hat Y^{M}=Y^{\Gamma^M(\tau_1^*,\ldots,\tau_j^*);\beta_1^*,\ldots,\beta_j^*}$ and using \eqref{ekv:rcINt1} of property \ref{vfass:rcINt}) we find that\\ $\sup_{s\in[0,T]}|Y_s^{\tau^*_1,\ldots,\tau^*_j;\beta^*_1,\ldots,\beta^*_j}-\hat Y^{M}_{s}|\to 0$ in probability, as $M\to\infty$. Hence, there is a subsequence $(M_k)_{k\geq 1}$ such that the limit taken over the subsequence is 0, $\Prob$-a.s. Furthermore, as the convergence is uniform the limit process is \cadlagp% For all $\tau^*_{j}\leq t\leq s $ we have
%\begin{align*}
%Y^{\tau^*_1,\ldots,\tau^*_j;\beta^*_1,\ldots,\beta^*_j}_t&=\mathop{\lim}_{k\to\infty}\hat Y^{M_k}_t\geq \mathop{\lim}_{k\to\infty}\E\Big[\hat Y^{M_k}_s\big|\mcF_t\Big]=E\Big[Y^{\tau^*_1,\ldots,\tau^*_j;\beta^*_1,\ldots,\beta^*_j}_s\big|\mcF_t\Big]
%\end{align*}
%$\Prob$-a.s.~where we have used the supermartingale property to reach the inequality and dominated convergence (see property~\ref{vfass:Ybnd})) to reach the last equality.

By right-continuity of the switching costs and $\Psi$ and \eqref{ekv:rcINt2} of property \ref{vfass:rcINt}) we have that $\E[\sup_{s\in[0,T]}|U_s-\hat U^{M_k}_{s}|^2]\to 0$ as $k\to\infty$, where for notational simplicity we abuse the notation in \eqref{ekv:dominated} and let
\begin{align*}
U:=\Big(&\ett_{[s=T]}\Psi(\tau^*_1,\ldots,\tau^*_j;\beta^*_1,\ldots,\beta^*_j)+\ett_{[s < T]}\max_{\beta\in\mcI^{-\beta^*_j}}\Big\{-c_{\beta^*_j,\beta}(s)
+Y^{\tau^*_1,\ldots,\tau^*_j,s;\beta^*_1,\ldots,\beta^*_j,\beta}_s\Big\}:\: \tau^*_{j}\leq s\leq T\Big).
\end{align*}
Hence, $(M_k)_{k\geq 0}$ has a subsequence $(\tilde M_k)_{k\geq 0}$ such that $\sup_{s\in[0,T]}|U_s-\hat U^{\tilde M_k}_{s}|\to 0$, $\Prob$-a.s.~as $k\to\infty$. This implies that $U$ is a \cadlag process which is of class [D] by Assumption~\ref{ass:prelim}.(\ref{ass:psiBND}) and property \ref{vfass:Ybnd}).

We thus have that $\hat U^{\tilde M_k}$ is a sequence of \cadlag processes of class [D] that converges pointwisely to the \cadlag process $U$ of class [D] and that $\hat Y^{\tilde M_k}$ is the Snell envelope of $\hat U^{\tilde M_k}$, for all $k\geq 0$. Then by Theorem~\ref{thm:Snell}.(\ref{Snell:lim}) we find that $\hat Y^{\tilde M_k}$ converges pointwisely to the Snell envelope Snell envelope of $U$. Hence, $\Big(Y^{\tau^*_1,\ldots,\tau^*_j;\beta^*_1,\ldots,\beta^*_j}_s:\: \tau^*_{j}\leq s\leq T\Big)$ is the Snell envelope of $U$.

To arrive at the second equality in \eqref{ekv:esssupATT} we note that the results we obtained in Step 1 %and quasi left-continuity of $\bbF$
implies that for any sequence $(\gamma_l)_{l\geq 0}\subset\mcT$ with $\gamma_l\nearrow\gamma\in\mcT$ we have $\lim_{l\to\infty}\E[\hat U^{M}_{\gamma_l}]\leq \E[\hat U^{M}_{\gamma}]$ for all $M\geq 1$. Now, for all $k\geq 0$ this gives
\begin{align*}
\lim_{l\to\infty}\E[U_{\gamma_l}]&\leq \lim_{l\to\infty}\E[\hat U^{\tilde M_k}_{\gamma_l}]+\lim_{l\to\infty}\E[|U_{\gamma_l}-\hat U^{\tilde M_k}_{\gamma_l}|]
\leq \E[U_{\gamma}]+2\E[\sup_{s\in[0,T]}|U_{s}-\hat U^{\tilde M_k}_{s}|],
\end{align*}
where the last term can be made arbitrarily small and we, thus, have that $\lim_{l\to\infty}\E[U_{\gamma_l}]\leq \E[U_{\gamma}]$ and by Theorem~\ref{thm:Snell}.(\ref{Snell:att}) we get \eqref{ekv:esssupATT}.
%\begin{align*}
%Y^{\tau^*_1,\ldots,\tau^*_j;\beta^*_1,\ldots,\beta^*_j}_{\tau^*_j} &=\E\Big[\ett_{[\tau^*_{j+1}=T]}\Psi(\tau^*_1,\ldots,\tau^*_j;\beta^*_1,\ldots,\beta^*_j)
%\\
%&\quad+\ett_{[\tau^*_{j+1} < T]}\left\{-c_{\beta^*_j,\beta^*_{j+1}}(\tau^*_{j+1}) + Y^{\tau^*_1,\ldots,\tau^*_{j+1};\beta^*_1,\ldots,\beta^*_{j+1}}_{\tau^*_{j+1}}\right\}\Big| \mcF_{\tau^*_j}\Big]
%\end{align*}
%$\Prob$-a.s.

By induction we get that for each $K\geq 0$,
\begin{align*}
Y_0&=\E\Big[\ett_{[N^*\leq K]}\Psi(\tau^*_1,\ldots,\tau^*_{N^*};\beta^*_1,\ldots,\beta^*_{N^*})-\sum_{j=1}^{K\wedge N^*}c_{\beta^*_{j-1},\beta^*_{j}}(\tau^*_j)
\\
&\qquad+\ett_{[N^*> K]}\{-c_{\beta^*_{K},\beta^*_{K+1}}(\tau^*_{K+1})+Y^{\tau^*_1,\ldots,\tau^*_{K+1};\beta^*_1,\ldots,\beta^*_{K+1}}_{\tau^*_{K+1}}\Big].
\end{align*}
Now, arguing as in the proof of Proposition~\ref{prop:finSTRAT} and using property \emph{\ref{vfass:Ybnd})} we find that $u^*\in\mcU^f$. Letting $K\to\infty$ and using dominated convergence we conclude that $Y_0=J(u^*)$.\\

\noindent {\bf Step 3} It remains to show that the strategy $u^*$ is optimal. %, \ie that for any strategy $u\in\mcU$ we have $J(u)\leq J(u^*)$.
To do this we pick any other strategy $\hat u:=(\hat\tau_1,\ldots,\hat\tau_{\hat N};\hat\beta_1,\ldots,\hat\beta_{\hat N})\in\mcU^f$. By the definition of $Y_0$ in~\eqref{ekv:Ydef} we have
\begin{align*}
Y_0 &\geq \E\Big[\ett_{[\hat\tau_1 \geq T]}\Psi_0 + \ett_{[\hat \tau_1 < T]}\max_{\beta\in \mcI^{-b_0}}\left\{-c_{b_0,\beta}(\hat\tau_1)+Y^{\hat \tau_1;\beta}_{\hat \tau_1}\right\}\Big]
\\
&\geq\E\Big[\ett_{[\hat\tau_1 \geq T]}\Psi_0 + \ett_{[\hat \tau_1 < T]}\left\{-c_{b_0,\hat\beta_1}(\hat\tau_1)+Y^{\hat \tau_1;\hat\beta_1}_{\hat \tau_1}\right\}\Big]
\end{align*}
but in the same way
\begin{align*}
Y^{\hat \tau_1,\hat \beta_1}_{\hat \tau_1}\geq\E\Big[\ett_{[\hat\tau_2 \geq T]}\Psi(\hat\tau_1,\hat\beta_1) + \ett_{[\hat \tau_2 < T]}\left\{-c_{\hat\beta_1,\hat\beta_2}(\hat\tau_2)+Y^{\hat \tau_1,\hat \tau_2;\hat\beta_1,\hat\beta_2}_{\hat \tau_1}\right\}\Big| \mcF_{\hat\tau_1}\Big],
\end{align*}
$\Prob$--a.s. By repeating this argument and using the dominated convergence theorem we find that $J(u^*)\geq J(\hat u)$ which proves that $u^*$ is in fact optimal. Repeating the above procedure with $(\vect,\vecb)\in\mcD^f$ as initial condition \eqref{ekv:Yuni} follows.\qed

\bigskip

%\begin{rem}\end{rem}
The main difference between the above proof and the proof of Theorem 1 in the original work by Djehiche, Hamad\'ene and Popier \cite{BollanMSwitch1} is that, due to the fact that the future reward at any time depends on the entire history of the control, we are forced consider a family of processes indexed by an uncountable set rather than a $q$-tuple for some finite positive $q$. Hence, we cannot simply write $Y^{\tau_1^*,\ldots,\tau_j^*;\beta_1^*,\ldots,\beta_j^*}$ as the sum of a finite number of Snell envelopes. To arrive at the above verification theorem we therefore impose the right-continuity constraint assumed in Assumption~\ref{ass:prelim}.\ref{ass:psiRCint}. This effectively allowed us to find the two sequences of processes that approach on the one hand the value process corresponding to the optimal control and on the other hand the dominated process, in $\mcS^2$.

%%%%%%%%%%%%%%%%%%%%%%%%%%%%%%%%%%%%%%%%%%%%%%%%%%%%%%%%%%%%%%%%%%%%%%%%%%%%%%%%%%%%%%%%%%%%%%%%%%%%%%%%%%%%%%%%%%%%%%%%%%%%%%%%%%%%%%%%%%%%%%%%%%
%%%%%%%%%%%%%%%%%%%%%%%%%%%%%%%%%%%%%%%%%%%%%%%%%%%%%%%%%%%%%%%%%%%%%%%%%%%%%%%%%%%%%%%%%%%%%%%%%%%%%%%%%%%%%%%%%%%%%%%%%%%%%%%%%%%%%%%%%%%%%%%%%%
%%%%%%%%%%%%%%%%%%%%%%%%%%%%%%%%%%%%%%%%%%%%%%%%%%%%%%%%%%%%%%%%%%%%%%%%%%%%%%%%%%%%%%%%%%%%%%%%%%%%%%%%%%%%%%%%%%%%%%%%%%%%%%%%%%%%%%%%%%%%%%%%%%

\section{Existence\label{sec:exist}}
Theorem~\ref{thm:vfc} presumes existence of the verification family $((Y^{\vect;\vecb}_s)_{0\leq s\leq T}: (\vect,\vecb)\in \bigSET^f)$. To obtain a satisfactory solution to Problem 1, we thus need to establish that a verification family exists. This is the topic of the present section. We will follow the standard existence proof which goes by applying a Picard iteration (see \cite{CarmLud,BollanMSwitch1,HamZhang}). We thus define a sequence $((Y^{\vect;\vecb,k}_s)_{0\leq s\leq T}: (\vect,\vecb)\in \mcD^f)_{k\geq 0}$ of families of processes as
\begin{align}
Y^{\vect;\vecb,0}_s:=\E\Big[\Psi(\vect;\vecb)\Big| \mcF_s\Big]\label{ekv:Y0def}
\end{align}
and
\begin{align}
Y^{\vect;\vecb,k}_s:=\esssup_{\tau \in \mcT_{s\vee t_n}} \E\Big[&\ett_{[\tau\geq T]}\Psi(\vect;\vecb)
+\ett_{[\tau < T]}\max_{\beta\in\mcI^{-b_n}}\left\{-c_{b_n,\beta}(\tau)+Y^{\vect,\tau;\vecb,\beta,k-1}_\tau\right\}\Big| \mcF_s\Big]\label{ekv:Ykdef}
\end{align}
for $k\geq 1$.

\bigskip

\begin{prop}\label{prop:Ykbnd}
The sequence $((Y^{\vect;\vecb,k}_s)_{0\leq s\leq T}: (\vect,\vecb)\in \mcD^f)_{k\geq 0}$ is uniformly bounded in the sense that there is a $K>0$ such that,
\begin{align*}
\sup_{u\in\mcU}\E\Big[\sup_{s\in[0,T]}|Y^{\tau_1,\ldots;\beta_1,\ldots,k}_s|^2\Big]\leq K,
\end{align*}
and for all $(\vect,\vecb)\in\mcD^f$ and $b\in\mcI^{-b_n}$, we have
\begin{align*}
\E\Big[\sup_{s\in[0,T]}|Y^{\vect,s\vee t_n;\vecb,b,k}_s|^2\Big]\leq K,
\end{align*}
for all $k\geq 0$.
\end{prop}

\noindent\emph{Proof.} By the definition of $Y^{\vect;\vecb,k}$ we have that for any $u\in\mcU^f$,
\begin{align*}
\E\Big[\Psi(\tau_1,\ldots;\beta_1,\ldots)\big|\mcF_s\Big]\leq Y^{\tau_1,\ldots;\beta_1,\ldots,k}_s&\leq \esssup_{\hat u\in\mcU}\E\Big[\Psi(\hat \tau_1,\ldots;\hat \beta_1,\ldots)\big|\mcF_s\Big].
\end{align*}
By Doob's maximal inequality we have that for any $\hat u:=(\hat \tau_1,\ldots;\hat \beta_1,\ldots)\in\mcU$
\begin{align*}
\E\Big[\sup_{s\in[0,T]}\E\Big[|\Psi(\hat \tau_1,\ldots;\hat \beta_1,\ldots)|\big|\mcF_s\Big]^2\Big]\leq C\E\Big[|\Psi(\hat \tau_1,\ldots;\hat \beta_1,\ldots)|^2\Big].
\end{align*}
Taking the supremum over all $\hat u\in\mcU$ on both sides and using that the right hand side is uniformly bounded by Assumption~\ref{ass:onPSIandC}.(\ref{ass:psiBND}.a) the first bound follows.\\

Concerning the second claim, note that
\begin{align*}
&\E\Big[\sup_{s\in[0,T]}|Y^{\vect,s\vee t_n;\vecb,b,k}_s|^2\Big]
\leq\sup_{u\in\mcU} \E\Big[\sup_{s\in[0,T]}\E[\sup_{r\in[t_n,T]}|\Psi(\vect,r,\tau_1\vee r,\ldots;\vecb,b,\beta_1,\ldots)|\big|\mcF_s]^2\Big].
\end{align*}
Now, arguing as above we find that
\begin{align*}
\E\Big[\sup_{s\in[0,T]}|Y^{\vect,s\vee t_n;\vecb,b,k}_s|^2\Big]
\leq C\sup_{u\in\mcU}\E\Big[\sup_{r\in[t_n,T]}|\Psi(\vect,r,\tau_1\vee r,\ldots;\vecb,b,\beta_1,\ldots)|^2\Big]
\end{align*}
where the right hand side is bounded by Assumption~\ref{ass:onPSIandC}.(\ref{ass:psiBND}.b).
\qed

\bigskip

\begin{prop}\label{prop:Yk}
The family of processes $((Y^{\vect;\vecb,k}_s)_{0\leq s\leq T}: (\vect,\vecb)\in \mcD^f)$ satisfies:
\begin{enumerate}[i)]
  \item\label{Yk:rcINt} For every $n\geq 1$ and every $(\eta,\vecb)\in \bar\mcT^n\times\bar\mcI^n$ and $b\in\mcI^{-b_n}$ we have
  \begin{equation*}
  \E\big[\sup_{s\in[0,T]}|Y^{\Gamma^l(\eta);\vecb,k}_s-Y^{\eta;\vecb,k}_s|^2\big]\to 0
  \end{equation*}
  and
  \begin{equation*}
  \E\big[\sup_{s\in[0,T]}|Y^{\Gamma^l(\eta),s\vee\Gamma^l(\eta_n);\vecb,b_n,k}_s-Y^{\eta,s\vee\eta_n;\vecb,b_n,k}_s|^2\big]\to 0,
  \end{equation*}
  as $l\to\infty$ uniformly in $k$.
  \item\label{Yk:qlc} For every $(\vect,\vecb)\in \mcD^f$ and every $b\in\I^{-b_n}$, the process $(Y^{\vect,s\vee t_n;\vecb,b,k}_s:0\leq s\leq T)$ is in $\Sqlc^2$ for $k=0,1,\ldots$
\end{enumerate}
\end{prop}

\noindent\emph{Proof.} The proof will follow by induction and we use (i') to denote the first statement without the uniformity.\\

For $k=0$, we have $Y^{\vect,\cdot\vee t_n;\vecb,b,0}_\cdot=V_\cdot^{\vect,\cdot\vee t_n;\vecb,b,\emptyset}\in\Sqlc^2$ by Assumption~\ref{ass:prelim}.(\ref{ass:rBNDYbnd}) and (\ref{Yk:rcINt}') follows from Assumption~\ref{ass:prelim}.(\ref{ass:psiRCint}). Now, assume that there is a $k'\geq 0$ such that (i') and (\ref{Yk:qlc}) holds for all $k\leq k'$. Applying a reasoning similar to that in the proof of Theorem~\ref{thm:vfc} we find that
\begin{align*}
Y^{\vect;\vecb,k'+1}_s=\esssup_{u \in \mcU^{k'+1}_{s\vee t_n}} V^{\vect;\vecb,u}_s.
\end{align*}
But then by Assumption~\ref{ass:prelim} we find that (\ref{Yk:rcINt}') and (\ref{Yk:qlc}) hold for $k'+1$. By induction (i') and (\ref{Yk:qlc}) hold for all $k\geq 0$.\\

It remains to show that (\ref{Yk:rcINt}) holds. By the above reasoning we find that, for each $k$ we have
\begin{align*}
&\E\big[\sup_{s\in[0,T]}|Y^{\Gamma^l(\eta);\vecb,k}_s-Y^{\eta;\vecb,k}_s|^2\big]
\\
&\leq \E\big[\sup_{s\in[0,T]}|(Y^{\Gamma^l(\eta);\vecb,k}_s-Y^{\eta;\vecb,k}_s)^+|^2\big] + \E\big[\sup_{s\in[0,T]}|(Y^{\eta;\vecb,k}_s - Y^{\Gamma^l(\eta);\vecb,k}_s)^+|^2\big]
\\
&\leq \sup_{u\in\mcU_{\Gamma^l(\eta_n)}}\E\big[\sup_{s\in[0,T]}|(V^{\Gamma^l(\eta);\vecb,u}_s-V^{\eta;\vecb,u}_s)^+|^2\big]
+ \sup_{u\in\mcU}\E\big[\sup_{s\in[0,T]}|(V^{\eta;\vecb,u}_s-V^{\Gamma^l(\eta);\vecb,\hat u^l}_s)^+|^2\big]
\end{align*}
where the right hand side of the last inequality does not depend on $k$ and tends to zero as $l\to\infty$ by Assumption~\ref{ass:prelim}.(\ref{ass:psiRCint}). The second statement in (\ref{Yk:rcINt}) follows by an identical argument.\qed

\begin{cor}
For each $k\geq 0$ and each $s\in[0,T]$ there is a $u^k=(\tau^k_1,\ldots,\tau^k_{N^k};$ $\beta^k_1,\ldots,\beta^k_{N^k})\in\mcU^k_{t_n\vee s}$, such that
\begin{align*}
Y^{\vect;\vecb,k}_s=\E\bigg[\Psi(\vect,\tau^k_1,\ldots,\tau^k_{N^k};\vecb,\beta^k_1,\ldots,\beta^k_{N^k}) -\sum_{j=1}^{N^k}c_{\beta^k_j,\beta^k_{j-1}}(\tau^k_j)\Big|\mcF_s\bigg],
\end{align*}
with $\beta^k_0=b_0$.
\end{cor}
\noindent\emph{Proof.} Follows from the definition of $Y^{\vect;\vecb,k}$ and Propositions~\ref{prop:Ykbnd} and \ref{prop:Yk} by applying the same argument as in the proof of the verification theorem (Theorem~\ref{thm:vfc}).\\

\bigskip

\begin{prop}\label{prop:Yklim}
For each $(\vect,\vecb)\in\mcD^f$, the limit $\bar Y^{\vect;\vecb}:=\lim_{k\to\infty}Y^{\vect;\vecb,k}$, exists as an increasing pointwise limit, $\Prob$-a.s. Furthermore, the process $\bar Y^{\vect,\cdot\vee t_n;\vecb,b}_\cdot$ is \cadlag for each $b\in\mcI^{-b_n}$.
\end{prop}

\noindent\emph{Proof.} Since $\mcU^k_t\subset \mcU^{k+1}_t$ we have that, $\Prob$-a.s.,
\begin{align*}
Y^{\vect;\vecb,k}_s \leq Y^{\vect;\vecb,k+1}_s\leq \esssup_{u\in\mcU}\E\Big[|\Psi(\tau_1,\ldots;\beta_1,\ldots)|\big|\mcF_s\Big],
\end{align*}
where the right hand side is bounded $\Prob$-a.s.~by Proposition~\ref{prop:Ykbnd}. Hence, the sequence $((Y^{\vect;\vecb,k}_s)_{0\leq s\leq T}: (\vect,\vecb)\in \bigSET)$ is increasing and $\Prob$-a.s.~bounded, thus, it converges $\Prob$-a.s.~for all $s\in [0,T]$.\\

Concerning the second claim, note that for $p\in(1,2)$, we have
\begin{align*}
\sup_{s\in[0,T]}Y^{\vect,s\vee t_n;\vecb,b,k}_s&\leq\sup_{s\in[0,T]}\sup_{r\in[0,T]}Y^{\vect,r\vee t_n;\vecb,b,k}_s
\\
&\leq \sup_{s\in[0,T]}\esssup_{u\in\mcU}\E[\sup_{r\in[t_n,T]}|\Psi(\vect,r,\tau_1\vee r,\ldots;\vecb,b,\beta_1,\ldots)|\big|\mcF_s]
\\
&\leq1+\sup_{s\in[0,T]}\esssup_{u\in\mcU}\E[\sup_{r\in[t_n,T]}|\Psi(\vect,r,\tau_1\vee r,\ldots;\vecb,b,\beta_1,\ldots)|^p\big|\mcF_s]=: K(\omega)
\end{align*}
for all $k\geq 0$ (where the inequalities hold $\Prob$-a.s.). Now, arguing as in the proof of Proposition~\ref{prop:Ykbnd} we have
\begin{align*}
&\E\big[\sup_{s\in[0,T]}\esssup_{u\in\mcU}\E[\sup_{r\in[t_n,T]}|\Psi(\vect,r,\tau_1\vee r,\ldots;\vecb,b,\beta_1,\ldots)|^p\big|\mcF_s]^{2/p}\big]
\\
&\leq C\sup_{u\in\mcU}\E\big[\sup_{r\in[t_n,T]}|\Psi(\vect,r,\tau_1\vee r,\ldots;\vecb,b,\beta_1,\ldots)|^2\big]<\infty.
\end{align*}
We thus conclude that there is a $\Prob$-null set $\mcN$ such that for each $\omega\in\Omega\setminus\mcN$ we have $K(\omega)<\infty$.

By the ``no-free-loop'' condition (Assumption~\ref{ass:onPSIandC}.(\ref{ass:oncNFL})) and the finiteness of $\mcI$ we get that for any control $(\tau_1,\ldots,\tau_N;\beta_1,\ldots,\beta_N)$,
\begin{align*}
\sum_{j=1}^{N}c_{\beta_j,\beta_{j-1}}(\tau_j)\geq \epsilon (N-m)/m,
\end{align*}
$\Prob$-a.s. For $\omega\in\Omega\setminus\mcN$ (in the remainder of the proof $\mcN$ denotes a generic $\Prob$-null set), we thus have
\begin{align*}
-K(\omega)\leq Y^{\vect,s\vee t_n;\vecb,b,k}_s(\omega)&\leq \E[\Psi(\vect,s\vee t_n,\tau^k_1,\ldots,\tau^k_{N^k};\vecb,b,\beta_1,\ldots,\beta^k_{N^k})
-\epsilon(N^k/m-1)|\mcF_s](\omega)
\\
&\leq K(\omega)+\epsilon-\epsilon/m\E[ N^k |\mcF_s](\omega),
\end{align*}
where $(\tau^k_1,\ldots,\tau^k_{N^k};\beta^k_1,\ldots,\beta^k_{N^k})\in\mcU_{s\vee t_n}$ is a control corresponding to $Y^{\vect,s\vee t_n;\vecb,b,k}_s$. This implies that for $k'>0$ we have,
\begin{align*}
\Prob[N^k>k' |\mcF_s](\omega)\leq (2K(\omega)m/\epsilon+m)/k'.
\end{align*}
Now, for all $0\leq k'\leq k$ we have,
\begin{align*}
&\breve Y^{\vect,s\vee t_n;\vecb,b,k,k'}_s:=\E\bigg[\Psi(\vect,s,\tau^k_1,\ldots,\tau^k_{N^k\wedge k'};\vecb,b,\beta^k_1,\ldots,\beta^k_{N^k\wedge k'})
- \sum_{j=1}^{N^k\wedge k'}c_{\beta^k_{j-1},\beta^k_{j}}(\tau^k_j)\Big|\mcF_s\bigg]
\\
&\leq Y^{\vect,s\vee t_n;\vecb,b,k'}_s\leq Y^{\vect,s\vee t_n;\vecb,b,k}_s,
\end{align*}
where we introduced the process $\breve Y^{\vecb,\vect,k,k'}$ corresponding to the truncation $(\tau^k_1,\ldots,\tau^k_{N^k\wedge k'};\beta^k_1,\ldots,\beta^k_{N^k\wedge k'})$ of the optimal control. As the truncation only affects the performance of the controller when $N^k>k'$ we have
\begin{align*}
Y^{\vect,s\vee t_n;\vecb,b,k}_s-\breve Y^{\vect,s\vee t_n;\vecb,b,k,k'}_s
&=\E\bigg[\ett_{[N^k>k']}\Big(\Psi(\vect,s\vee t_n,\tau^k_1,\ldots,\tau^k_{N^k};\vecb,b,\beta^k_1,\ldots,\beta^k_{N^k}) - \sum_{j=1}^{N^k}c_{\beta^k_{j-1},\beta^k_{j}}(\tau^k_j)
\\
&\quad-\Psi(\vect,s\vee t_n,\tau^k_1,\ldots,\tau^k_{N^k\wedge k'};\vecb,b,\beta^k_1,\ldots,\beta^k_{N^k\wedge k'}) + \sum_{j=1}^{N^k\wedge k'}c_{\beta^k_{j-1},\beta^k_{j}}(\tau^k_j)\Big)\Big|\mcF_s\bigg]
\\
&\leq\E\bigg[\ett_{[N^k>k']}\Big(\Psi(\vect,s\vee t_n,\tau^k_1,\ldots,\tau^k_{N^k};\vecb,b,\beta^k_1,\ldots,\beta^k_{N^k})
\\
&\quad-\Psi(\vect,s\vee t_n,\tau^k_1,\ldots,\tau^k_{N^k\wedge k'};\vecb,b,\beta^k_1,\ldots,\beta^k_{N^k\wedge k'}) \Big)\Big|\mcF_s\bigg].
\end{align*}
Applying H\"older's inequality we get that for $\omega\in\Omega\setminus\mcN$,
\begin{align*}
&Y^{\vect,s\vee t_n;\vecb,b,k}_s(\omega)-\breve Y^{\vect,s\vee t_n;\vecb,b,k,k'}_s(\omega)
\\
&\leq 2\E[\ett_{[N^k>k']}|\mcF_s]^{1/q}
\esssup_{u\in\mcU}\E[\sup_{r\in[t_n,T]}|\Psi(\vect,r,\tau_1\vee r,\ldots;\vecb,b,\beta_1,\ldots)|^p\big|\mcF_s]^{1/p}(\omega)
\\
&\leq 2((K(\omega)m/\epsilon+m)/k')^{1/q}(K(\omega))^{1/p},
\end{align*}
with $\frac{1}{p}+\frac{1}{q}=1$, there is thus a constant $C=C(\omega)$ such that
\begin{align*}
Y^{\vect,s\vee t_n;\vecb,b,k}_s(\omega)-Y^{\vect,s\vee t_n;\vecb,b,k'}_s(\omega)\leq C(k')^{-1/q},
\end{align*}
for all $s\in[0,T]$. We conclude that for all $\omega\in\Omega\setminus\mcN$, the sequence\\ $(Y^{\vect,\cdot\vee t_n;\vecb,b,k}_\cdot(\omega))_{k\geq 0}$ is a sequence of \cadlag functions that converges uniformly which implies that the limit is a \cadlag function.\qed

\bigskip

\begin{prop}
The family $((\bar Y^{\vect;\vecb}_s)_{0\leq s\leq T}:(\vect,\vecb)\in\mcD^f)$ is a verification family.
\end{prop}

\noindent\emph{Proof.} As $\bar Y^{\vect;\vecb}$ is the pointwise limit of an increasing sequence of \cadlag supermartingales it is a \cadlag supermartingale (see p.~86 in \cite{DelMeyer2}). We treat each remaining property in the definition of a verification family separately:\\

\noindent\emph{a)} Applying the convergence result to the right hand side of \eqref{ekv:Ykdef} and using the fact that, by Proposition~\ref{prop:Yklim},
\begin{align*}
\ett_{[s\geq T]}\Psi(\vect;\vecb)+\ett_{[s < T]}\max_{\beta\in\mcI^{-b_n}}\left\{-c_{b_n,\beta}(s)+\bar Y^{\vect,s\vee t_n;\vecb,\beta}_s\right\}
\end{align*}
is a \cadlag process, (iv) of Theorem \ref{thm:Snell} gives
\begin{align*}
\bar Y^{\vect;\vecb}_s:=\esssup_{\tau \in \mcT_{s}} \E\Big[&\ett_{[\tau\geq T]}\Psi(\vect;\vecb)+\ett_{[\tau < T]}\max_{\beta\in\mcI^{-b_n}}\left\{-c_{b_n,\beta}(\tau)+\bar Y^{\vect,\tau;\vecb,\beta}_\tau\right\}\Big| \mcF_s\Big].
\end{align*}

\bigskip

\noindent\emph{b)} Uniform boundedness was shown in Proposition~\ref{prop:Ykbnd}.

\bigskip

\noindent\emph{c)} We have
\begin{align*}
\lim_{l\to\infty}\E\big[\sup_{s\in[0,T]}|\bar Y^{\Gamma^l(\eta);\vecb}_s-\bar Y^{\eta;\vecb}_s|^2\big]&=\lim_{l\to\infty}\E\big[\sup_{s\in[0,T]}\lim_{k\to\infty}|\bar Y^{\Gamma^l(\eta);\vecb,k}_s-\bar Y^{\eta;\vecb,k}_s|^2\big]
\\
&\leq \lim_{l\to\infty}\lim_{k\to\infty}\E\big[\sup_{s\in[0,T]}|\bar Y^{\Gamma^l(\eta);\vecb,k}_s-\bar Y^{\eta;\vecb,k}_s|^2\big]
\\
&= \lim_{k\to\infty}\lim_{l\to\infty}\E\big[\sup_{s\in[0,T]}|\bar Y^{\Gamma^l(\eta);\vecb,k}_s-\bar Y^{\eta;\vecb,k}_s|^2\big]
\\
&=0
\end{align*}
where taking limits is interchangeable due to the uniform convergence property shown in Proposition~\ref{prop:Yk}.(\ref{Yk:rcINt}). The second statement in c), that is equation \eqref{ekv:rcINt2}, follows by an identical argument.

\bigskip

\noindent\emph{d)} We know from Proposition~\ref{prop:Yklim} that $\bar Y^{\vect,\cdot\vee t_n;\vecb,b}_\cdot$ is \cadlag and by Proposition~\ref{prop:Ykbnd} it follows that $\bar Y^{\vect,\cdot\vee t_n;\vecb,b}_\cdot\in\mcS^2$. It remains to show that $\bar Y^{\vect,\cdot\vee t_n;\vecb,b}_\cdot$ is quasi-left continuous. Using the notation from the proof of Proposition~\ref{prop:Yklim} we have for $k\geq 0$,
\begin{align*}
&|\bar Y^{\vect,\gamma_j(\omega)\vee t_n;\vecb,b}_{\gamma_j(\omega)}(\omega)-\bar Y^{\vect,\gamma(\omega)\vee t_n;\vecb,b}_{\gamma(\omega)}(\omega)|
\leq
|Y^{\vect,\gamma_j(\omega)\vee t_n;\vecb,b,k}_{\gamma_j(\omega)}(\omega) -Y^{\vect,\gamma(\omega)\vee t_n;\vecb,b,k}_{\gamma(\omega)}(\omega)|+2C(\omega)k^{-1/q},
\end{align*}
for all $\omega\in\Omega\setminus\mcN$ with $\Prob(\mcN)=0$. By Proposition~\ref{prop:Yk}.(\ref{Yk:qlc}) the first part tends to zero $\Prob$-a.s.~as $j\to \infty$. Since $k$ was arbitrary and $C$ is $\Prob$-a.s.~bounded the desired result follows. This finishes the proof.\qed

\section{Application to SDDEs with controlled volatility\label{sec:sdde}}
We now move to the case of impulse control of SDDEs. However, we start by formalizing the hydro-power production problem proposed as a motivating example in the introduction.

\subsection{Continuous time hydro-power planning}
The increasing competitiveness of electricity markets calls for new operational standards in electric power production facilities. It has previously been acknowledged that optimal switching can be useful in deriving production schedules that maximize the revenue from electricity production~\cite{CarmLud,BollanMSwitch1,KharrSC}. Here we will extend the applicability of optimal switching by introducing a new example, the coordinated operation of hydropower plants  interconnected by hydrological coupling.

We consider the situation where a central operator controls the output of two hydropower stations located in the same river (but note that the model is easily extended to consider an entire system of power stations).\\

We assume that Plant $i$, for $i=1,2$, has:
\begin{itemize}
  \item A reservoir containing a volume $Z^i_t$ m$^3$ of water at time $t$.
  \item A stochastic inflow $V^i_t$ m$^3$/s to the reservoir that is modeled by a jump diffusion process.
  \item $\kappa_i$ turbines that can be either ``in operation'', producing $p_i(Z^i_t)$ MW by releasing $\alpha_i$ m$^3$/s of water through the turbine or ``idle''.
\end{itemize}
We assume that the power plants are hydrologically connected in such a way that the water that passes through Plant 1 will reach the reservoir of Plant 2 after $\delta\geq 0$ seconds.\\

We assume that we control the number of turbines in operation in each of the two plants. We thus let $\mcI:=\{0,1,\ldots,\kappa_1\}\times\{0,1,\ldots,\kappa_2\}$. The dynamics of the involved processes is then given by
\begin{align*}
dV_t&=a(t,V_t)dt+\sigma(t,V_t)dW_t+\int_{\R^2\setminus\{0\}}\gamma(t,V_{t-},z)\Gamma(dt,dz)
\\
dZ^1_t&=(V^1_t-\alpha_1\xi^1_t)dt
\\
dZ^2_t&=(V^2_t-\alpha_2\xi^2_t+\alpha_1\xi^1_{t-\delta})dt
\\
(V_0,Z_0)&=(v_0,z_0)\in\R_+^4
\end{align*}
and an appropriate reward functional is
\begin{align*}
J(u):=\E\Big[\int_0^TR_t(\xi^1_t p_1(Z^1_t)+\xi^2_t p_2(Z^2_t))dt+q(Z^1_T,Z^2_T)\Big],
\end{align*}
where $R_t$ is the (stochastic) electricity price at time $t$ and $q:\R_ +^2\to\R$ is the value of water (per m$^3$) stored in the reservoirs at the end of the operation period\footnote{Note that we expect the water in Reservoir 1 to have a higher value as it can be used in both plants.}.

\subsection{A general SDDE model}
Motivated by the above example we assume that $\bbF$ is the completed filtration generated by an $d$-dimensional Brownian motion $W$ and an $\dimX$-dimensional, independent, finite activity, Poisson random measure $\Gamma$ with intensity measure $\nu(ds; dz) = ds \times \mu(dz)$, where $\mu$ is the L\'evy measure on $\R^\dimX$ of $\Gamma$ and  $\tilde \Gamma(ds; dz) := (\Gamma - \nu)(ds; dz)$ is called the compensated jump martingale
random measure of $\Gamma$. For $u\in\mcU$, we let $X^{u,0}$ solve
\begin{align}
dX^{u,0}_t&=a(t,X^{u,0}_t,X^{u,0}_{t-\delta})dt+\sigma(t,X^{u,0}_t,X^{u,0}_{t-\delta})dW_t\nonumber
\\
&\quad+\int_{\R^\dimX\setminus\{0\}} \gamma(t,X^{u,0}_{t-},X^{u,0}_{t-\delta},z)\tilde \Gamma(dt,dz),\quad{\rm for\, all}\,t\in(0,T],\label{ekv:SDDE0_1}
\\
X^{u,0}_{s}&=\chi(s),\quad s\in [-\delta,0],\label{ekv:SDDE0_2}
\end{align}
where $\delta>0$ is a constant and $\chi:[-\delta,0]\to\R^\dimX$ is a deterministic \cadlag function with $\sup_{s\in [-\delta,0]}|\chi(s)|\leq C$, and define recursively
\begin{align}
dX^{u,j}_t&=a(t,X^{u,j}_t,X^{u,j}_{t-\delta})dt+\sigma(t,X^{u,j}_t,X^{u,j}_{t-\delta})dW_t\nonumber
\\
&\quad+\int_{\R^\dimX\setminus\{0\}} \gamma(t,X^{u,j}_{t-},X^{u,j}_{t-\delta},z)\tilde \Gamma(dt,dz),\quad{\rm for\, all}\,t\in(\tau_{j},T],\label{ekv:SDDE1}
\\
X^{u,j}_{\tau_j}&=h_{\beta_{j-1},\beta_j}(\tau_j,X^{u,j-1}_{\tau_j})\label{ekv:SDDE2}
\\
X^{u,j}_{s}&=X^{u,j-1}_{s},\quad s\in [-\delta,\tau_j).\label{ekv:SDDE3}
\end{align}
Finally we let $X^u:=\lim_{j\to\infty}X^{u,j}$ be our controlled process\footnote{Whenever it exists, we refer to the limit process $X^u$ as a solution to the SDDE \eqref{ekv:SDDE1}-\eqref{ekv:SDDE3}}.

\begin{rem}
Note that by letting $\chi_1\equiv b_0$ and taking $[h_{\beta_{j-1},\beta_j}]_1(t,x)=\beta_j$ and letting the first rows of $a$, $\sigma$ and $\gamma$ equal zeros we get $[X]_1=\xi^u$ which implies that the control enters all terms in the SDDE for $X^u$.
\end{rem}

We consider the situation when the functional $J$ is given by
\begin{align*}
J(u):=\E\bigg[\int_0^T f(t,X^u_t)dt+g(X^u_T)-\sum_{j=1}^N c_{\beta_{j-1},\beta_j}(\tau_j)\bigg].
\end{align*}

We assume that the parameters of the SDDE satisfies the following conditions:
\begin{ass}\label{ass:onSDDE}
\begin{enumerate}[i)]
\item The functions $a:[0,T]\times\R^\dimX\times\R^\dimX\to\R^\dimX$ and $\sigma:[0,T]\times\R^\dimX\times\R^\dimX\to\R^\dimX\times\R^\dimX$ are continuous in $t$ and satisfy
\begin{align*}
|a(t,x,y)-a(t,x',y')|+|\sigma(t,x,y)-\sigma(t,x',y')|\leq C(|x-x'|+|y-y'|)
\end{align*}
for all $(x,x',y,y')\in \R^{4\dimX}$.
\item There is a $\rho(z)$, with $\int \rho^{4q}(z)\mu(dz)< \infty$ such that $\gamma:[0,T]\times\R^\dimX\times\R^\dimX\times\R^\dimX\to\R^\dimX$ satisfies
\begin{align*}
|\gamma(t,x,y,z)-\gamma(t,x',y',z)|&\leq \rho(z)(|x-x'|+|y-y'|),
\\
|\gamma(t,x,y,z)|&\leq \rho(z)(1+|x|+|y|).
\end{align*}
\item\label{ass:onSDDE.h} For all $(t,x)\in[0,T]\times\R^\dimX$ and all $(b,b')\in\bar\mcI^2$, the map $h_{b,b'}:[0,T]\times\R^\dimX\to\R^\dimX$ satisfies
\begin{align*}
|h_{b,b'}(t,x)|\leq C\vee |x|.
\end{align*}
Furthermore,
\begin{align*}
|h_{b,b'}(t,x)-h_{b,b'}(t',x')|\leq |x-x'|+C|t-t'|
\end{align*}
for all $(x,x')\in \R^{2\dimX}$ and $(t,t')\in[0,T]^2$.
\end{enumerate}
\end{ass}

\begin{rem}
Note in particular that since $a$ and $\sigma$ are continuous in $t$, $a(\cdot,0,0)$ and $\sigma(\cdot,0,0)$ are uniformly bounded and Lipschitz continuity implies that
\begin{align}
&|a(t,x,y)|^{4q}+|\sigma(t,x,y)|^{4q}+\int_{\R^\dimX\setminus\{0\}} |\gamma(t,x,y,z)|^{4q}\mu(dz)
\leq C(1+|x|^{4q}+|y|^{4q}).\label{ekv:coeffBND}
\end{align}
\end{rem}

We have the following result:
\begin{prop}\label{prop:SDDExnu}
Under Assumption~\ref{ass:onSDDE} the SDDE \eqref{ekv:SDDE1}-\eqref{ekv:SDDE3} admits a unique solution for each $u\in\mcU$. Furthermore, the solution has moments of order $4q$, \ie $\sup_{u\in\mcU}\E\big[\sup_{t\in[0,T]}|X^u_t|^{4q}\big]<\infty$.
\end{prop}

\noindent\emph{Proof.} We first note that existence of a unique solution to the SDDE follows by repeated use of Theorem 3.2 in~\cite{AgramOksen} (where existence of a unique solution to a more general controlled SDDE is shown). It remains to show that the moment estimate holds. We have $X^{u,j}=X^{u,j-1}$ on $[-\delta,\tau_{j})$ and
\begin{align*}
X^{u,j}_t&=h_{\beta_{j-1},\beta_{j}}(\tau_{j},X^{u,j-1}_{\tau_{j}})+\int_{\tau_{j}}^ta(s,X^{u,j}_s,X^{u,j}_{s-\delta})ds
\\
&\quad + \int_{\tau_{j}}^t\sigma(t,X^{u,j}_s,X^{u,j}_{s-\delta})dW_s+\int_{\tau_{j}}^t \int_{\R^\dimX\setminus\{0\}} \gamma(s,X^{u,j}_{s-},X^{u,j}_{s-\delta},z)\tilde \Gamma(ds,dz)
\end{align*}
on $[\tau_{j},T]$. By Assumption~\ref{ass:onSDDE}.(\ref{ass:onSDDE.h}) we get, for $t\in [\tau_{j},T]$, using integration by parts, that
\begin{align*}
|X^{u,j}_t|^2&= |X^{u,j}_{\tau_{j}}|^2+2\int_{\tau_{j}+}^t X^{u,j}_{s-} dX^{u,j}_s+\int_{\tau_{j}+}^t d[X^{u,j},X^{u,j}]_s
\\
&\leq C\vee|X^{u,j-1}_{\tau_{j}}|^2+2\int_{\tau_{j}+}^t X^{u,j}_{s-} dX^{u,j}_s+\int_{\tau_{j}+}^t d[X^{u,j},X^{u,j}]_s
\\
&\leq C\vee|X^{u,j-1}_{\tau_{j-1}}|^2+2\int_{\tau_{j-1}+}^{\tau_{j}} X^{u,j-1}_{s-} dX^{u,j-1}_s+\int_{\tau_{j-1}+}^{\tau_{j}}d[X^{u,j-1},X^{u,j-1}]_s
\\
&\quad+2\int_{\tau_{j}+}^t X^{u,j}_{s-} dX^{u,j}_s+\int_{\tau_{j}+}^t d[X^{u,j},X^{u,j}]_s.
\end{align*}
By repeated application we find that
\begin{align*}
|X^{u,j}_t|^2&\leq C\vee|X^{u,0}_{0}|^2+\sum_{i=0}^{j-1} \{2\int_{\tau_{i}+}^{\tau_{i+1}} X^{u,i}_{s-} dX^{u,i}_s+\int_{\tau_{i}+}^{\tau_{i+1}} d[X^{u,i},X^{u,i}]_s\}
\\
&\quad+2\int_{\tau_{j}+}^t X^{u,j}_{s-} dX^{u,j}_s+\int_{\tau_{j}+}^t d[X^{u,j},X^{u,j}]_s
\\
&\leq C+\sum_{j=0}^{j-1} \big\{2\int_{\tau_{i}+}^{\tau_{i+1}} X^{u,i}_{s-} dX^{u,i}_s+\int_{\tau_{i}+}^{\tau_{i+1}} d[X^{u,i},X^{u,i}]_s\big\}
\\
&\quad+2\int_{\tau_{j}+}^t X^{u,j}_{s-} dX^{u,j}_s+\int_{\tau_{j}+}^t d[X^{u,j},X^{u,j}]_s,
\end{align*}
with $\tau_0:=0$. Now, since $X^{u,i}$ and $X^{u,j}$ coincide on $[0,\tau_{i+1\wedge j+1})$ we have
\begin{align*}
\sum_{i=0}^{j-1} \int_{\tau_{i}+}^{\tau_{i+1}} X^{u,i}_{s-} dX^{u,i}_s+\int_{\tau_{j}+}^t X^{u,j}_{s-} dX^{u,j}_s
&=\int_{0}^{t}X_{s}^{u,j}a(s,X^{u,j}_s,X^{u,j}_{s-\delta})ds+\int_{0}^{t}X_{s}^{u,j}\sigma(s,X^{u,j}_s,X^{u,j}_{s-\delta})dW_s\\
&\quad+\int_{0}^{t}\int_{\R^\dimX\setminus\{0\}}X_{s-}^{u,j} \gamma(s,X^{u,j}_{s-},X^{u,j}_{s-\delta},z)\tilde \Gamma(ds,dz)
\end{align*}
and
\begin{align*}
&\E\Big[\sum_{i=0}^{j-1} \int_{\tau_{i}+}^{\tau_{i+1}} d[X^{u,i},X^{u,i}]_s + \int_{\tau_{j}+}^t d[X^{u,j},X^{u,j}]_s\Big]
\\
&=\E\Big[\int_{0}^{t}(|a(s,X^{u,j}_s,X^{u,j}_{s-\delta})|^2+\int_{\R^\dimX\setminus\{0\}}|\gamma(s,X^{u,j}_{s-},X^{u,j}_{s-\delta},z)|^2\mu(dz))ds\Big].
\end{align*}
Finally, using the Burkholder-Davis-Gundy inequality in combination with \eqref{ekv:coeffBND} we get
\begin{align*}
\E\Big[\sup_{s\in[0,t]}|X^{u,j}_s|^{4q}\Big]&\leq C + C\int_{0}^t\E\Big[\sup_{r\in[0,s]}|X^{u,j}_r|^{4q}\Big]ds,
\end{align*}
where the constant $C$ does not depend on $j$ and it follows by Gr\"onwall's lemma that $\E\Big[\sup_{t\in[0,T]}|X^{u,j}_t|^{4q}\Big]$ is bounded uniformly in $j$. Now, the result follows since $\tau_j\to T$, $\Prob$-a.s.,~as $j\to\infty$.\qed

\bigskip

For each $(\vect,\vecb)\in\mcD^f$ and each $u\in\mcU$ we let
\begin{align*}
X^{\vect;\vecb,u}:=X^{t_1,\ldots,t_n,t_n\vee\tau_1,\ldots,t_n\vee\tau_N;b_1,\ldots,b_n,\beta_1,\ldots,\beta_N}
\end{align*}
and% similarly we let
\begin{align*}
X^{\vect;\vecb,u,j}:=X^{t_1,\ldots,t_n,t_n\vee\tau_1,\ldots,t_n\vee\tau_{N};b_1,\ldots,b_n,\beta_1,\ldots,\beta_{N},j}.
\end{align*}

\begin{prop}\label{prop:SDDEbndB}
For all $(\vect,\vecb)\in\mcD^f$ we have
\begin{align*}
\sup_{u\in\mcU}\E\big[\sup_{s\in[0,T]}\sup_{t\in[t_n,T]}|X^{\vect,t;\vecb,b,u}_s|^{4q}\big]<\infty.
\end{align*}
\end{prop}

\noindent\emph{Proof.} %Let for $r\in [0,T]$, let $(Z^{r}_t:0\leq t\leq T)$ be defined as $Z^{r}_t:=X^{}$ on $[0,r)$ and
For $t\in[t_n,T]$ we have, for $s\geq t$,
\begin{align*}
X^{\vect,t;\vecb,b}_s&=h_{b_n,b}(t,X^{\vect;\vecb}_t)+\int_{t}^s a(r,X^{\vect,t;\vecb,b}_r,X^{\vect,t;\vecb,b}_{r-\delta})dr
\\
&\quad + \int_{t}^s\sigma(r,X^{\vect,t;\vecb,b}_r,X^{\vect,t;\vecb,b}_{r-\delta})dW_r
+\int_{t}^s \int_{\R^\dimX\setminus\{0\}} \gamma(r,X^{\vect,t;\vecb,b}_{r-},X^{\vect,t;\vecb,b}_{r-\delta},z)\tilde \Gamma(dr,dz).
\end{align*}
Arguing as in the proof of Proposition~\ref{prop:SDDExnu} we find that for $s\in [\tau_{j},T]$,
\begin{align*}
\sup_{t\in[t_n,T]}|X^{\vect,t;\vecb,b,u,n+1+j}_s|^2
&\leq C\vee\sup_{t\in[t_n,T]}|X^{\vect;\vecb}_{t}|^2
+\sup_{t\in[t_n,T]}\Big\{\sum_{i=0}^{j-1} \big\{2\int_{t\vee\tau_{i}+ }^{\tau_{i+1}} X^{\vect,t;\vecb,b,u,n+1+i}_{r-} dX^{\vect,t;\vecb,b,u,n+1+i}_r
\\
&\quad+\int_{\tau_{i}+}^{\tau_{i+1}} d[X^{\vect,t;\vecb,b,u,n+1+i},X^{\vect,t;\vecb,b,u,n+1+i}]_r\big\}
\\
&\quad+2\int_{t\vee\tau_{j}+}^s X^{\vect,t;\vecb,b,u,n+1+j}_{r-} dX^{\vect,t;\vecb,b,u,n+1+j}_{r}
\\
&\quad+\int_{\tau_{j}+}^s d[X^{\vect,t;\vecb,b,u,n+1+j},X^{\vect,t;\vecb,b,u,n+1+j}]_r\Big\}.
\end{align*}
We thus find that, for each $u\in\mcU$,
\begin{align*}
&\E\Big[\sup_{s\in[0,T]}\sup_{t\in[t_n,T]}|X^{\vect,t;\vecb,b,u}_s|^{4q}\Big]
\leq C +C\E\Big[\sup_{s\in[0,t]}|X^{\vect;\vecb}_s|^{4q}\Big] + C\int_{0}^t\E\Big[\sup_{s\in[0,T]}\sup_{t\in[t_n,T]}|X^{\vect,t;\vecb,b,u}_s|^{4q}\Big]ds
\end{align*}
and the assertion again follows by applying Gr\"onwall's lemma and using Proposition~\ref{prop:SDDExnu}.\qed

To illustrate that switching does not diverge solutions we have the following useful lemma:
\begin{lem}\label{lem:SDDEcont}
For $\gamma\in\mcT$ and each $u\in\mcU_{\gamma}$, let $(^k\! Z^u)_{k\geq 0}$ and $X^u$ be processes in $\mcS^{4q}$ (with $\E[\sup_{s\in[0,\gamma]}|^k Z^u|^{4q}]$ uniformly bounded) that solve the SDDE \eqref{ekv:SDDE1}-\eqref{ekv:SDDE3} on $(\gamma,T]$ with control $u$ and such that
\begin{align}
\E\Big[\int_0^{\gamma}|X_s^{u}-^k\!\!Z^{u}_s|^{4}ds+|X_{\gamma}^{u,0}-^k\!\!Z^{u,0}_\gamma|^{4}\Big]\to 0,
\end{align}
as $k\to\infty$. Then,
\begin{align}\label{ekv:lemCont1}
\lim_{k\to\infty}\sup_{u\in\mcU_\gamma}\E\Big[\sup_{s\in[\gamma,T]}|X_s^{u}-^k\!\!Z^{u}_s|^{2}\Big]\to 0
\end{align}
and for all $b\in\mcI^{-b_0}$ we have
\begin{align}\label{ekv:lemCont2}
\lim_{k\to\infty}\sup_{u\in\mcU_\gamma}\E\Big[\sup_{t\in[\gamma,T]}\sup_{s\in[\gamma,T]}|X_s^{t,b,u}-^k\!\!Z^{t,b,u}_s|^{2}\Big]\to 0.
\end{align}
\end{lem}

\noindent\emph{Proof.} %For $t\geq\gamma$. We have
%\begin{align*}
%X_t^{u,0}&= X_t^{u,0}+\int_0^t a(s,X_s^{u,0},X_{s-\delta}^{u,0})ds+\int_0^t\sigma(s,X_s^{u,0},X_{s-\delta}^{u,0})dB_s
%\\
%&\quad+\int_0^t\int_{\R^\dimX\setminus{0}} \gamma(s,X_s^{u,0},X_{s-\delta}^{u,0})\tilde \Gamma(ds,dz).
%\end{align*}
%We get
%\begin{align*}
%X_t^{u,0}-^k\!\!Z^{u,0}_t&= X_{\gamma}^{u,0}-^k\!\!Z^{u,0}_\gamma+\int_0^t(a(s,X_s^{u,0},X_{s-\delta}^{u,0})- a(s,X_s^{\tilde \chi,u,0},^k\!\!Z_{s-\delta}^{u,0}))ds
%\\
%&+\int_0^t(\sigma(s,X_s^{u,0},X_{s-\delta}^{u,0})-\sigma(s,X_s^{\tilde \chi,u,0},^k\!\!Z_{s-\delta}^{u,0}))dB_s
%\\
%&+\int_0^t\int_{\R^\dimX\setminus{0}} (\gamma(s,X_{s-}^{u,0},X_{s-\delta}^{u,0})-\gamma(s,^k\!\!Z_{s-}^{u,0},^k\!\!Z_{s-\delta}^{u,0}))\tilde \Gamma(ds,dz)
%\end{align*}
%We thus get (using \eg Thm 66, p. 339 in~\cite{Protter} and Lipschitz continuity)
%\begin{align*}
%\E\Big[\sup_{s\in[0,t]} |X_s^{u,0}-X^{\tilde \chi,u,0}_s|^2\Big]&\leq C(|\chi(0)-\tilde\chi(0)|^2+\E\Big[\int_0^{t}(|X_s^{u,0}-X_s^{\tilde \chi,u,0}|^2+|X_{s-\delta}^{u,0}-{}^k\!\!Z_{s-\delta}^{u,0}|^2)ds\big]).
%\end{align*}
By the contraction property of $h_{.,.}$ we have that $|X_{\tau_j}^{u,j}-^k\!\!Z_{\tau_j}^{u,j}|<|X_{\tau_j}^{u,j-1}-^k\!\!Z_{\tau_j}^{u,j-1}|$. Using integration by parts we get, for $t\in[\tau_j,T]$,
\begin{align*}
|X_{t}^{u,j}-^k\!\!Z_{t}^{u,j}|^2
& = |X_{\tau_j}^{u,j}-^k\!\!Z_{\tau_j}^{u,j}|^2
+2\int_{\tau_j+}^t(X_{s-}^{u,j}-^k\!\!Z_{s-}^{u,j})(dX_{s}^{u,j}-d^k\!Z_{s}^{u,j})
\\
&\quad+\int_{\tau_j+}^t d[X^{u,j}-^k\!\!Z^{u,j},X^{u,j}-^k\!\!Z^{u,j}]_s
\\
&\leq |X_{\tau_{j-1}}^{u,j-1}-^k\!\!Z_{\tau_{j-1}}^{u,j-1}|^2+2\int_{\tau_{j-1}}^{\tau_{j}}(X_{s-}^{u,j-1} -^k\!\!Z_{s-}^{u,j-1})(dX_{s}^{u,j-1}-d^k\!Z_{s}^{u,j-1})
\\
&\quad+2\int_{\tau_j+}^t(X_{s-}^{u,j}-^k\!\!Z_{s-}^{u,j})(dX_{s}^{u,j}-d^k\!Z_{s}^{u,j})
+\int_{\tau_{j-1}+}^{\tau_{j}} d[X^{u,j-1}-^k\!\!Z^{u,j-1},X^{u,j-1}-^k\!\!Z^{u,j-1}]_s
\\
&\quad+\int_{\tau_j+}^t d[X^{u,j}-^k\!\!Z^{u,j},X^{u,j}-^k\!\!Z^{u,j}]_s.
\end{align*}
Repeated application implies that
\begin{align*}
|X_{t}^{u}-^k\!\!Z_{t}^{u}|^2& \leq |X_{\gamma}^{u,0}-^k\!\!Z_{\gamma}^{u,0}|^2+2\sum_{j=0}^{\infty}\int_{\tau_j+}^{\tau_{j+1}\wedge t}(X_{s-}^{u,j}-^k\!\!Z_{s-}^{u,j})(dX_{s}^{u,j}-d^k\!Z_{s}^{u,j})
\\
&\quad+\sum_{j=0}^{\infty}\int_{\tau_j+}^{\tau_{j+1}\wedge t} d[X^{u,j}-^k\!\!Z^{u,j},X^{u,j}-^k\!\!Z^{u,j}]_s.
\end{align*}
Now, for $s\in (\tau_j,T]$ we have
\begin{align*}
dX_{s}^{u,j}-d^k\!Z_{s}^{u,j}&=(a(s,X^{u,j}_s,X^{u,j}_{s-\delta})-a(s,^k\!\!Z^{u,j}_s,^k\!\!Z^{u,j}_{s-\delta}))ds
\\
&\quad+(\sigma(s,X^{u,j}_s,X^{u,j}_{s-\delta})-\sigma(s,^k\!\!Z^{u,j}_s,^k\!\!Z^{u,j}_{s-\delta}))dW_s\\
&\quad+\int_{\R^\dimX\setminus\{0\}} (\gamma(s,X^{u,j}_{s-},X^{u,j}_{s-\delta},z)-\gamma(s,^k\!\!Z^{u,j}_{s-},^k\!\!Z^{u,j}_{s-\delta},s))\tilde \Gamma(ds,dz).
\end{align*}
Using Lipschitz continuity of $a,\sigma$ and $\gamma$ and the Burkholder-Davis-Gundy inequality we get
\begin{align*}
\E\Big[\sup_{s\in[\gamma,t]}|X_{s}^{u}-^k\!\!Z_{s}^{u}|^4\Big]& \leq C\E\Big[|X_{\gamma}^{u,0}-^k\!\!Z_{\gamma}^{u,0}|^4+\int_0^{\gamma}|X_{s}^{u}-^k\!\!Z_{s}^{u}|^4ds\Big]
+ C\int_{\gamma}^{t}\E\Big[\sup_{r\in[\gamma,s]}|X_{r}^{u}-^k\!\!Z_{r}^{u}|^4\Big]ds,
\end{align*}
where the constant $C$ does not depend on the  control $u$, and by Gr\"onwall's inequality we have
\begin{align*}
\E\Big[\sup_{s\in[\gamma,t]}|X_{s}^{u}-^k\!\!Z_{s}^{u}|^4\Big]& \leq C\E\Big[|X_{\gamma}^{u,0}-^k\!\!Z_{\gamma}^{u,0}|^4+\int_0^{\gamma}|X_{s}^{u}-^k\!\!Z_{s}^{u}|^4ds\Big].
\end{align*}
Now, applying Jensen's inequality gives \eqref{ekv:lemCont1}. Furthermore, we have
\begin{align*}
\sup_{r\in[0,T]}|X_{t}^{r,b,u}-^k\!\!Z_{t}^{r,b,u}|^2 &\leq \sup_{r\in[0,T]}|X_{r}^{u,0}-^k\!\!Z_{r}^{u,0}|^2
\\
&\quad+2\sup_{r\in[0,T]}\Big\{\sum_{j=0}^{\infty}\int_{\tau_j+\vee r}^{\tau_{j+1}\wedge t}(X_{s-}^{r,b,u,j}-^k\!\!Z_{s-}^{r,b,u,j})(dX_{s}^{r,b,u,j}-d^k\!Z_{s}^{r,b,u,j})
\\
&\quad+\sum_{j=0}^{\infty}\int_{\tau_j+}^{\tau_{j+1}\wedge t} d[X^{r,b,u,j}-^k\!\!Z^{r,b,u,j},X^{r,b,u,j}-^k\!\!Z^{r,b,u,j}]_s\Big\}.
\end{align*}
and \eqref{ekv:lemCont2} follows by an identical argument.\qed

\bigskip

We add the following assumptions on the components of the cost functional and the functions $h$.
\begin{ass}\label{ass:onfgc}
\begin{enumerate}[(i)]
  \item The functions $f:[0,T]\times \R^\dimX\to\R$ and $g:\R^\dimX\to\R$ are both locally Lipschitz in $x$. Furthermore, there are constants $q> 1$ and $K>0$ such that
      \begin{equation*}
      |f(t,x)|+|g(x)|\leq K(1+|x|^q)
      \end{equation*}
      for all $(t,x)\in [0,T]\times\R^\dimX$.
  \item For all $b\in\mcI$ we have
  \begin{equation*}
    g(x)>\max_{b'\in\mcI^{-b}}g(h_{b,b'}(T,x))-c_{b,b'}(T),
  \end{equation*}
  for all $x\in\R^\dimX$.
%  \item\label{ass:onc2} The switching costs $(c_{b,b'})_{b,b'\in \mcI}\in (\mcS_{qc}^2)^{m\times m}$ are such that
%  \begin{enumerate}[(a)]
%    \item $\inf_{t\in [0,T]} c_{b,b'}(t)\geq 0$,\quad for all $(b,b')\in \mcI\times\mcI$
%    \item $c_{b_{1},b_2}(t_1)+c_{b_2,b_3}(t_2)+\ldots+c_{b_{n-1},b_n}(t_n)+c_{b_{k},b_1}(t)\geq\delta>0$,\quad for all $(t_1,\ldots,t_n,b_1,\ldots,b_n)\in\mcD^f$ and all $t\in [t_n,T]$.
%  \end{enumerate}
  \item There is a constant $\kappa>0$ such that for any sequence $(b_1,\ldots,b_{j})\in\bar\mcI^j$ with $j>\kappa$ there is a subsequence $1 = \iota_1<\cdots<\iota_{j'}=j$ with $j'\leq\kappa$ and $(b_{\iota_1},\ldots,b_{\iota_{j'}})\in\bar\mcI^{j'}$ for which
\begin{align*}
&h_{b_1,b_2}(t,h_{b_2,b_3}(t,\cdots h_{b_{j-1},b_{j}}(t,x)\cdots)
=h_{b_{\iota_1},b_{\iota_2}}(t,h_{b_{\iota_2},b_{\iota_3}}(t,\cdots h_{b_{\iota_{j'-1}},b_{\iota_{j'}}}(t,x)\cdots)).
\end{align*}
\end{enumerate}
\end{ass}
It is straightforward to see that with the above assumptions the $\Psi$ defined by
\begin{equation*}
\Psi(\vect;\vecb):=\int_0^T f(t,X^{\vect;\vecb}_t)dt+g(X^{\vect;\vecb}_T)
\end{equation*}
satisfies Assumption~\ref{ass:onPSIandC}. The remainder of this section is devoted to showing that $\Psi$ also satisfies Assumption~\ref{ass:prelim}, guaranteeing the existence of an optimal control to the problem of maximizing $J$.

\begin{prop}
For each $n\geq 1$ and each $(\eta,\vecb)\in\bar\mcT^n\times\bar\mcI^n$ and $b\in\mcI^{-b_n}$ there is a map $(\mcU\to\mcU:u\to \hat u^l)_{l\geq 1}$ such that
\begin{align}
\lim_{l\to\infty}\sup_{u\in\mcU}\E\Big[\sup_{s\in[0,T]}|(V^{\eta;\vecb,u}_s-V^{\Gamma^l(\eta);\vecb,\hat u^l}_s)^+|^2\Big]=0\label{ekv:VrcInETA1}
\end{align}
and% a map $(\mcU\to\mcU:u\to \hat u^l)_{l\geq 1}$ such that
\begin{align}
\lim_{l\to\infty}\sup_{u\in\mcU}\E\Big[\sup_{s\in[0,T]}|(V^{\eta,s\vee \eta_n;\vecb,b,u}_s-V^{\Gamma^l(\eta),s\vee \Gamma^l(\eta_n);\vecb,b,\hat u^l}_s)^+|^2\Big]=0.\label{ekv:VrcInETA2}
\end{align}
Furthermore, we have
\begin{align}
\lim_{l\to\infty}\sup_{u\in\mcU_{\Gamma^l(\eta_n)}}\E\Big[\sup_{s\in[0,T]}|(V^{\Gamma^l(\eta);\vecb,u}_s-V^{\eta;\vecb,u}_s)^+|^2\Big]=0 \label{ekv:VrcInETA3}
\end{align}
and
\begin{align}
\lim_{l\to\infty}\sup_{u\in\mcU_{\Gamma^l(\eta_n)}}\E\Big[\sup_{s\in[0,T]}|(V^{\Gamma^l(\eta),s\vee \Gamma^l(\eta_n);\vecb,b,u}_s-V^{\eta,s\vee \eta_n;\vecb,b,u}_s)^+|^2\Big]=0. \label{ekv:VrcInETA4}
\end{align}
\end{prop}

\noindent\emph{Proof.} To simplify notation we let $(\zeta_i)_{1\leq i\leq n}$ denote %the non-decreasing sequence of $\bbF$-stopping times
$\Gamma^l(\eta)$ and let $X$ and $Z$ (resp.~$X^j$ and $Z^j$) denote $X^{\eta;\vecb,u}_t$ resp.~$X^{\Gamma^l(\eta);\vecb,\hat u^l}$ (resp.~$X^{\eta;\vecb,u,j}$ and $X^{\Gamma^l(\eta);\vecb,\hat u^l,j}$). Furthermore, we let $U^*_t:=\sup_{s\in[0,t]}|U_s|$ be the running maximum of the process $|U|$.

We have:\\

\emph{i)} $X_t= Z_t$, for all $t\in[0,\eta_1)$, $\Prob$-a.s.\\

\emph{ii)} On $[\eta_1,\zeta_1)$ we have $|X_t-Z_t|\leq (X)_T^*+(Z)_T^*$.\\

\emph{iii)} If $\eta_j\leq\zeta_1$, then $\zeta_j=\zeta_{j-1}=\cdots=\zeta_1$.\\

Letting $M_1:=\max\{j\geq 1: \eta_j\leq \zeta_1\}$ we get
\begin{align*}
X^{M_1}_{\zeta_{M_1}}-Z^{M_1}_{\zeta_{M_1}}&=X^{{M_1}}_{\zeta_{M_1}}+(h_{b_{M_1-1},b_{M_1}}(\eta_{M_1},X^{M_1-1}_{\eta_{M_1}}) -X^{{M_1}}_{\eta_{M_1}})
-h_{b_{M_1-1},b_{M_1}}(\zeta_{M_1},Z^{M_1-1}_{\zeta_{M_1}}).
\end{align*}
Hence,
\begin{align*}
|X^{M_1}_{\zeta_{M_1}}-Z^{M_1}_{\zeta_{M_1}}|&\leq |X^{{M_1}}_{\zeta_{M_1}}-X^{{M_1}}_{\eta_{M_1}}|+ C|\eta_{M_1}-\zeta_{M_1}| + |X^{M_1-1}_{\eta_{M_1}} - Z^{M_1-1}_{\zeta_{M_1}}|
\\
&\leq  C2^{-l} +|X^{{M_1}}_{\zeta_{M_1}}-X^{{M_1}}_{\eta_{M_1}}|+ |X^{{M_1-1}}_{\zeta_{M_1}}-X^{{M_1-1}}_{\eta_{M_1}}|
 + |X^{M_1-1}_{\zeta_{M_1}} - Z^{M_1-1}_{\zeta_{M_1}}|.
\end{align*}
But $X^{0}_{\zeta_1}=Z^{0}_{\zeta_1}$ and by induction it follows that
\begin{align*}
|X^{M_1}_{\zeta_{M_1}}-Z^{M_1}_{\zeta_{M_1}}|&\leq M_1C2^{-l}+\sum_{j=1}^{M_1}(|X^{{j}}_{\zeta_{j}}-X^{j}_{\eta_{j}}|+ |X^{{j-1}}_{\zeta_{j}}-X^{{j-1}}_{\eta_{j}}|).
\end{align*}

If we iteratively define $M_i:=\max\{j > M_{i-1} : \eta_j\leq \zeta_{M_{i-1}+1}\}$, for $i=1,\ldots n_M$ with $M_{n_M}= n$ and $M_0:=0$. Then we get, in the same manner,
\begin{align*}
|X^{M_i}_{\zeta_{M_i}}-Z^{M_i}_{\zeta_{M_i}}|&\leq (M_i-M_{i-1}) C2^{-l}+\!\!\!\sum_{j=M_{i-1}+1}^{M_i}(|X^{{j}}_{\zeta_{j}}-X^{j}_{\eta_{j}}|+ |X^{{j-1}}_{\zeta_{j}}-X^{{j-1}}_{\eta_{j}}|)
\\
&\quad + |X^{M_{i-1}}_{\zeta_{M_i}} - Z^{M_{i-1}}_{\zeta_{M_i}}|.
\end{align*}

Now on $[{\zeta_{M_i}},T]$ we have
\begin{align*}
X_t^{M_i}-Z^{M_i}_t&= X^{M_i}_{\zeta_{M_i}}-Z^{M_i}_{\zeta_{M_i}}+\int_{\zeta_{M_i}}^t(a(s,X^{M_i}_s,X^{M_i}_{s-\delta})-a(s,Z^{M_i}_s,Z^{M_i}_{s-\delta}))ds
\\
&\quad+\int_{\zeta_{M_i}}^t(\sigma(s,X^{M_i}_s,X_{s-\delta}^{M_i})-\sigma(s,Z_s^{M_i},Z_{s-\delta}^{M_i}))dB_s
\\
&\quad+\int_{\zeta_{M_i}}^t \int_{\R^\dimX\setminus\{0\}}(\gamma(s,X^{M_i}_{s-},X_{s-\delta}^{M_i})-\gamma(s,Z_{s-}^{M_i},Z_{s-\delta}^{M_i}))\tilde \Gamma(ds,dz).
\end{align*}
Put together we find that for $t\in[\zeta_{M_i},T]$ we have
\begin{align*}
|X_{t}^{M_i}-Z^{M_i}_t|&\leq (M_i-M_{i-1}) C2^{-l}+\!\!\!\sum_{j=M_{i-1}+1}^{M_i}(|X^{{j}}_{\zeta_{j}}-X^{j}_{\eta_{j}}|+ |X^{{j-1}}_{\zeta_{j}}-X^{{j-1}}_{\eta_{j}}|)
\\
&\quad + |X^{M_i-1}_{\zeta_{M_i}} - Z^{M_i-1}_{\zeta_{M_i}}|+\int_{\zeta_{M_i}}^t|a(s,X^{M_i}_s,X^{M_i}_{s-\delta})-a(s,Z^{M_i}_s,Z^{M_i}_{s-\delta})|ds
\\
&\quad+|\int_{\zeta_{M_i}}^t(\sigma(s,X^{M_i}_s,X_{s-\delta}^{M_i})-\sigma(s,Z_s^{M_i},Z_{s-\delta}^{M_i}))dB_s
\\
&\quad+\int_{\zeta_{M_i}}^t \int_{\R^\dimX\setminus\{0\}}(\gamma(s,X^{M_i}_{s-},X_{s-\delta}^{M_i})-\gamma(s,Z_{s-}^{M_i},Z_{s-\delta}^{M_i}))\tilde \Gamma(ds,dz)|.
\end{align*}
Applying Thm 66, p. 339 in~\cite{Protter} and Lipschitz continuity iteratively gives
\begin{align*}
\E\Big[\sup_{s\in[\zeta_{M_i},t]}|X_{s}^{M_i}-Z^{M_i}_s|^4\Big]&\leq C2^{-l}+C\E\Big[\sum_{j=1}^{M_i}(|X^{{j}}_{\zeta_{j}}-X^{j}_{\eta_{j}}|^4
\\
&\quad + |X^{{j-1}}_{\zeta_{j}}-X^{{j-1}}_{\eta_{j}}|^4)+ \int_{0}^t (|X^{M_i}_s-Z^{M_i}_s|^4+|X^{M_i}_{s-\delta}-Z_{s-\delta}^{M_i}|^4)ds\Big].
\end{align*}
By Gr\"onwall's inequality and point ii) above we find that
\begin{align}\nonumber
\E\Big[\sup_{t\in[\zeta_{M_i},T]}|X_{t}^{M_i}-Z^{M_i}_t|^4\Big]&\leq C2^{-l}(1+(X^*_T)^4+(Z^*_T)^4)
\\
&\quad+C\sum_{j=1}^{M_i}\E\big[|X^{{j}}_{\zeta_{j}}-X^{j}_{\eta_{j}}|^4+ |X^{{j-1}}_{\zeta_{j}}-X^{{j-1}}_{\eta_{j}}|^4\big].\label{ekv:asLEM}
\end{align}
Moving on we consider the possibility of interventions in the period $[\eta_n,\zeta_n)$. Let $N':=\max\{j\geq 0: \tau_j<\zeta_n\}$ and note that if $N'> \kappa$, then there is a subsequence $(\iota_j)_{j=1}^{\kappa'}$ with $1\leq \iota_1 < \cdots <\iota_{\kappa'}=N'$ with $\kappa'\leq\kappa$ and $(b_n,\beta_{\iota_1}, \ldots, \beta_{\iota_{\kappa'}})\in \bar\mcI^{\kappa'+1}$ such that, for all $(t,x)\in[0,T]\times \R^\dimX$,
\begin{align*}
h_{b_n,\beta_1}\circ\cdots\circ h_{\beta_{N'-1},\beta_{N'}}(t,x)=h_{b_n,\beta_{\iota_1}}\circ\cdots\circ h_{\beta_{\iota_{\kappa'-1}},\beta_{\iota_{\kappa'}}}(t,x).
\end{align*}
We then let\footnote{For $k\geq 1$ we denote by $\bold 1_k$ the vector of $k$ ones.} $\hat u^l=(\hat\tau_1,\ldots,\hat\tau_{\hat N};\hat\beta_1,\ldots,\hat\beta_{\hat N}):=$\\ $(\zeta_{n}{\bold 1}_{\kappa'},\tau_{N'+1},\ldots,\tau_N;\beta_{\iota_1},\ldots,\beta_{\iota_{\kappa'}},\beta_{N'+1},\ldots,\beta_N)$. Arguing as above, we find that
\begin{align}
|X_{\zeta_n}-Z_{\zeta_n}|&\leq N' C2^{-l}+\sum_{j=1}^{N'}(|X^{{n+j}}_{\zeta_n}-X^{n+j}_{\tau_{j}}|+ |X^{{n+j-1}}_{\zeta_n}-X^{{n+j-1}}_{\tau_{j}}|)
 +|X^{n}_{\zeta_n}-Z^{n}_{\zeta_n}|.\label{ekv:diff@n}
\end{align}
We now turn to the total revenue and let
\begin{equation*}
\Lambda:=\sum_{j=1}^{\hat N}c_{\hat \beta_{j-1},\hat \beta_j}(\hat \tau_j)-\sum_{j=1}^{N}c_{\beta_{j-1},\beta_j}(\tau_j).
\end{equation*}
By right continuity of the switching costs, we find that
\begin{align}
\lim_{l\to\infty}\Lambda\leq \bigg(\frac{\kappa}{2}-\frac{N'-m}{m}\bigg)\rho,\label{ekv:Cs}
\end{align}
$\Prob$-a.s. The difference in revenue can then be written
\begin{align*}
V^{\eta;\vecb,u}_t-V^{\zeta;\vecb,\hat u^l}_t&= \E\Big[\int_0^T (f(s,X_s)-f(s,Z_s))ds+g(X_T)-g(Z_T)+\Lambda\big| \mcF_t\Big].
\end{align*}
By local Lipschitz continuity of $f$ and $g$ we get that, for each $K>0$ there is a $C> 0$ such that $|f(t,x)-f(t,x')|\leq C|x-x'|$ and $|g(x)-g(x')|\leq C|x-x'|$ on $|x|+|x'|\leq K$. This gives us the relation
\begin{align*}
(V^{\eta;\vecb,u}_t-V^{\zeta;\vecb,\hat u^l}_t)^+
&\leq \E\Big[(\int_0^TC|X_s-Z_s|ds+C|X_T-Z_T|+\Lambda)^+\big| \mcF_t\Big]
\\
&\qquad+C\E[\ett_{[X_T^*+Z^*_T > K]}(1+(X_T^*)^q+(Z^*_T)^q )|\mcF_t]
\\
&\leq \E\Big[\ett_{A}(\int_0^TC|X_s-Z_s|ds+C|X_T-Z_T|+\Lambda^+)\big| \mcF_t\Big]
\\
&\qquad+C\E[\ett_{[X_T^*+Z^*_T > K]}(1+(X_T^*)^q+(Z^*_T)^q )|\mcF_t],
\end{align*}
where $A:=\{\omega\in\Omega:\int_0^TC|X_s-Z_s|^2ds+C|X_T-Z_T|^2>-\Lambda\}$. Doob's maximal inequality then gives that
\begin{align*}
\E\Big[\sup_{t\in[0,T]}((V^{\eta;\vecb,u}_t-V^{\zeta;\vecb,\hat u^l}_t)^+)^2\Big]
&\leq C\E\Big[\ett_{A}(\int_0^T|X_s-Z_s|^2ds+|X_T-Z_T|^2+(\Lambda^+)^2)\Big]
\\
&\quad+C\E[\ett_{[X_T^*+Z^*_T > K]}(1+(X_T^*)^{2q}+(Z^*_T)^{2q} )]
\\
&\leq C\E\Big[\ett_{A}(\int_0^T|X_s-Z_s|^2ds+|X_T-Z_T|^2+(\Lambda^+)^2)\Big]
\\
&\quad+C\Prob[X_T^*+Z^*_T > K]^{1/2},
\end{align*}
where we have used H\"older's inequality and the moment estimate in Proposition~\ref{prop:SDDExnu} to arrive at the last inequality.
For any $M>0$ we thus have
\begin{align}\nonumber
\E\Big[\sup_{t\in[0,T]}((V^{\eta;\vecb,u}_t-V^{\zeta;\vecb,\hat u^l}_t)^+)^2\Big]
&\leq C\E\Big[\ett_{[N'\leq M]}(\int_0^T|X_s-Z_s|^2ds+|X_T-Z_T|^2)\Big]\nonumber
\\
&\quad +C\E\Big[\ett_{[N'>M]}\ett_A((X_T^*)^2+(Z^*_T)^2)\Big]\nonumber
\\
&\quad+C\E\big[(\Lambda^+)^2\big]+C\Prob[X_T^*+Z^*_T > K]^{1/2},\label{ekv:bnd}
\end{align}
Concerning the first term, we have that $\ett_{[N'\leq M]}|X_s-Z_s|\leq |\tilde X_s-\tilde Z_s|$, where $\tilde X=X$ and $\tilde Z=Z$ on $[N'\leq M]$. On $[N'> M]$ we let $\tilde X:=X^{\eta;\vecb,\tilde u}$ with
\begin{align*}
\tilde u:=\left\{\begin{array}{cl} (\tau_1,\ldots,\tau_{M},\zeta_n,\tau_{N'+1},\ldots,\tau_N;\beta_1,\ldots,\beta_{M},\beta_{N'},\ldots,\beta_N), & {\rm if}\: \beta_{M}\neq\beta_{N'},
\\
(\tau_1,\ldots,\tau_{M},\tau_{N'+1},\ldots,\tau_N;\beta_1,\ldots,\beta_{M},\beta_{N'+1},\ldots,\beta_N), & {\rm if}\: \beta_{M}=\beta_{N'}.
\end{array}\right.
\end{align*}
and $\tilde Z:=X^{\eta;\vecb,\tilde u^l}$ where $\tilde u^l$ is obtained from $\tilde u$ as $\hat u^l$ was obtained from $u$. Now, we proceed as above and get for each $M\geq\kappa$, that
\begin{align*}
|\tilde X_{\zeta_n}-\tilde Z_{\zeta_n}|&\leq M C2^{-l}+\sum_{j=1}^{N'\wedge M}(|X^{{n+j}}_{\zeta_n}-X^{n+j}_{\tau_{j}}|+ |X^{{n+j-1}}_{\zeta_n}-X^{{n+j-1}}_{\tau_{j}}|)
 +|X^{n}_{\zeta_n}-Z^{n}_{\zeta_n}|.
\end{align*}
By \eqref{ekv:asLEM} and \eqref{ekv:lemCont1} of Lemma~\ref{lem:SDDEcont} we then find that for each $M\geq\kappa$, the first term on the right hand side in \eqref{ekv:bnd} goes to 0 as $l\to\infty$. Concerning the second term we have, again by H\"older's inequality and Proposition~\ref{prop:SDDExnu}, that
\begin{align*}
\E\Big[\ett_{[N'>M]}\ett_A((X_T^*)^2+(Z^*_T)^2)\Big]\leq C\Prob[[N'>M]\cap A]^{1/2}.
\end{align*}
Now, $A\subset \{\omega:C (X_T^*+Z^*_T)>-\Lambda\}$, where $C>0$ does not depend on $l$. For $l$ sufficiently large we thus see, by \eqref{ekv:Cs} and Chebyshev's inequality, that the probability on the right hand side can be made arbitrarily small by choosing $M$ sufficiently large. For the third term we note that
\begin{align*}
\E\big[(\Lambda^+)^2\big]\leq \kappa^2\sum_{(b,b')\in\bar\mcI^2}\E\big[\sup_{s\in[\eta_n,\zeta_n]}|c_{b,b'}(\zeta_n)-c_{b,b'}(s)|^2\big],
\end{align*}
where the right hand side goes to 0 as $l\to\infty$ by right-continuity of the switching costs. Finally, the last term of \eqref{ekv:bnd} can be made arbitrarily small by choosing $K$ large.

Concerning the second claim we note that with $X=X^{\eta,s\vee \eta_n,\vecb,b,u}$ and $Z=X^{\Gamma^l(\eta),s\vee \Gamma^l(\eta_n),\vecb,b,u}$ the relation in \eqref{ekv:diff@n} is replaced by
\begin{align*}
|X_{\zeta_n}-Z_{\zeta_n}|&\leq (N'+1) C2^{-l}+\sup_{r\in[\eta_1,\zeta_1]}\sum_{j=1}^{N'+1}(|X^{{n+j}}_{\zeta_n}-X^{n+j}_{r}|
\\ &\quad + |X^{{n+j-1}}_{\zeta_n}-X^{{n+j-1}}_{r}|)+|X^{n}_{\zeta_n}-Z^{n}_{\zeta_n}|.
\end{align*}
Hence, appealing to \eqref{ekv:lemCont2} of Lemma~\ref{lem:SDDEcont}, right-continuity and the result in Proposition~\ref{prop:SDDEbndB} the first second and last terms in the equivalent to \eqref{ekv:bnd} tends to 0 as $l\to\infty$ and \eqref{ekv:VrcInETA3} follows.

The last two statements given in equations \eqref{ekv:VrcInETA3}-\eqref{ekv:VrcInETA4} follow by a similar reasoning while noting that in this case $N'=0$ which implies that $\Lambda=0$, $\Prob$-a.s.\qed

\begin{lem}\label{lem:xCADLAG}
For all $(\vect,\vecb)\in\mcD^f$ and $k\geq 0$ we have
\begin{align*}
\sup_{u\in\mcU^k}\E\big[\sup_{s\in[t',T]}|X^{\vect,t';\vecb,b,u}_s-X^{\vect,t;\vecb,b,u}_s|\big|\mcF_{t'}\big]\to 0,
\end{align*}
$\Prob$-a.s.~as $t'\searrow t$.
\end{lem}

\noindent\emph{Proof.} Starting with $k=0$ we note that for $t'\geq t$ we have
\begin{align*}
X^{\vect,t;\vecb,b}_{t'}&=h_{b_n,b}(t,X^{\vect;\vecb}_{t})+X^{\vect,t;\vecb,b}_{t'}-X^{\vect,t;\vecb,b}_{t}%\int_{t}^{t'} a(r,X^{\vect,t;\vecb,b}_r,X^{\vect,t;\vecb,b}_{r-\delta})dr + \int_{t}^{t'}\sigma(t,X^{\vect,t;\vecb,b}_r,X^{\vect,t;\vecb,b}_{r-\delta})dW_r
%\\
%&\quad + \int_{t}^{t'}\int_{\R^\dimX\setminus\{0\}} \gamma(r,X^{\vect,t;\vecb,b}_{r-},X^{\vect,t;\vecb,b}_{r-\delta},z)\tilde \Gamma(dr,dz)).
\end{align*}
which gives
\begin{align*}
|X^{\vect,t';\vecb,b}_{t'}-X^{\vect,t;\vecb,b}_{t'}|&\leq C|t'-t|+|X^{\vect;\vecb}_{t'}-X^{\vect;\vecb}_{t}|+|X^{\vect,t;\vecb,b}_{t'}-X^{\vect,t;\vecb,b}_{t}|.
\end{align*}
For $k> 0$ and $u\in\mcU^k_t$ we have, for $i\leq k$
\begin{align*}
X^{\vect,t;\vecb,b,u,n+i+1}_{t'}&=\ett_{[\tau_i\leq t']}\{h_{\beta_{i-1},\beta_i}(\tau_i,X^{\vect,t;\vecb,b,u,n+i}_{\tau_i})+X^{\vect,t;\vecb,b,u,n+i+1}_{t'} \\
&\quad -X^{\vect,t;\vecb,b,u,n+i+1}_{\tau_i}\}+ \ett_{[\tau_i> t']}X^{\vect,t;\vecb,b,u,n+i}_{t'}
\end{align*}
and
\begin{align*}
X^{\vect,t';\vecb,b,u,n+i+1}_{t'}&=\ett_{[\tau_i\leq t']}h_{\beta_{i-1},\beta_i}(t',X^{\vect,t';\vecb,b,u,n+i}_{t'}) + \ett_{[\tau_i> t']}X^{\vect,t';\vecb,b,u,n+i}_{t'}.
\end{align*}
which gives
\begin{align*}
&|X^{\vect,t';\vecb,b,u,n+i+1}_{t'}-X^{\vect,t;\vecb,b,u,n+i+1}_{t'}|
\\
&\leq\ett_{[\tau_i\leq t']}\{C|t'-\tau_i|+|X^{\vect,t;\vecb,b,u,n+i}_{t'}-X^{\vect,t';\vecb,b,u,n+i}_{t'}|
\\
&\quad+|X^{\vect,t;\vecb,b,u,n+i}_{t'} - X^{\vect,t;\vecb,b,u,n+i}_{\tau_i}| + |X^{\vect,t;\vecb,b,u,n+i+1}_{t'} - X^{\vect,t;\vecb,b,u,n+i+1}_{\tau_i}|\} \\
&\quad+ \ett_{[\tau_i> t']}|X^{\vect,t;\vecb,b,u,n+i}_{t'}-X^{\vect,t';\vecb,b,u,n+i}_{t'}|.
%\int_{\tau_k}^s a(r,X^{\vect,t;\vecb,b,u}_r,X^{\vect,t;\vecb,b,u}_{r-\delta})dr + \int_{\tau_k}^s\sigma(t,X^{\vect,t;\vecb,b,u}_r,X^{\vect,t;\vecb,b,u}_{r-\delta})dW_r
%\\
%&\quad +\int_{\tau_k}^s \int_{\R^\dimX\setminus\{0\}} \gamma(r,X^{\vect,t;\vecb,b}_{r-},X^{\vect,t;\vecb,b}_{r-\delta},z)\tilde \Gamma(dr,dz)).
\end{align*}
Repeated application renders
\begin{align*}
|X^{\vect,t';\vecb,b,u}_{t'}-X^{\vect,t;\vecb,b,u}_{t'}|
&\leq C(k+1)|t'-t|+\sum_{i=1}^{k}\ett_{[\tau_i\leq t']}\{|X^{\vect,t;\vecb,b,u,n+i}_{t'} - X^{\vect,t;\vecb,b,u,n+i}_{\tau_i}|
\\
&+ |X^{\vect,t;\vecb,b,u,n+i+1}_{t'} - X^{\vect,t;\vecb,b,u,n+i+1}_{\tau_i}|\}+|X^{\vect;\vecb}_{t'}-X^{\vect;\vecb}_{t}|+|X^{\vect,t;\vecb,b}_{t'}-X^{\vect,t;\vecb,b}_{t}|.
\end{align*}
Furthermore, we have
\begin{align*}
\int_{0}^{t'}|X^{\vect,t';\vecb,b,u}_{s}-X^{\vect,t;\vecb,b,u}_{s}|^4ds&\leq |t'-t|((X^{\vect,t';\vecb,b,u})^*_T+(X^{\vect,t;\vecb,b,u})^*_T)^4,
\end{align*}
where the right hand side tends to zero $\Prob$-a.s. as $t'\searrow t$~by $\Prob$-a.s.~boundedness of\\ $\sup_{u\in\mcU}\sup_{r\in[t_n,T]}|(X^{\vect,r;\vecb,b,u})^*_T|^{4}$. Arguing as in the proof of Lemma~\ref{lem:SDDEcont} we find that
\begin{align*}
&\E\Big[\sup_{s\in[t',T]}|X^{\vect,t';\vecb,b,u}_{s}-X^{\vect,t;\vecb,b,u}_{s}|^4\big|\mcF_{t'}\Big] \leq C(|X^{\vect,t';\vecb,b,u}_{t'}-X^{\vect,t;\vecb,b,u}_{t'}|^4+\int_0^{t'}|X^{\vect,t';\vecb,b,u}_{s}-X^{\vect,t;\vecb,b,u}_{s}|^4ds),
\end{align*}
and the assertion follows by right continuity of $X$.\qed

\begin{lem}\label{lem:xLCE}
For all $(\vect,\vecb)\in\mcD^f$ and all $b\in \mcI^{-b_n}$ we have whenever $\gamma_j\nearrow\gamma\in\mcT_{t_n}$, with $(\gamma_j)_{j\geq 0}\subset\mcT_{t_n}$, that
\begin{align*}
\lim_{j\to\infty}\sup_{u\in\mcU^k_{\gamma_j}} \E\big[\sup_{s\in [\gamma,T]}|X_s^{\vect,\gamma_j;\vecb,b,u}-X_s^{\vect,\gamma;\vecb,b,u}|^2\big]=0,
\end{align*}
for all $0\leq k<\infty$.
\end{lem}

\noindent\emph{Proof.} Arguing as in the proof of the previous lemma we find that%, for $u\in\mcU_{\gamma_j}^k$, and $i\leq k$
%\begin{align*}
%|X^{\vect,\gamma_j;\vecb,b,u}_\gamma-X^{\vect,\gamma;\vecb,b,u}_\gamma|&\leq\ett_{[\tau_i\leq t']}\{C|\gamma-\gamma_j|+|X^{\vect,\gamma_j;\vecb,b,u,n+i}_{\gamma} - X^{\vect,\gamma_j;\vecb,b,u,n+i}_{\tau_i}|
%\\
%&\quad + |X^{\vect,\gamma_j;\vecb,b,u,n+i+1}_{\gamma} - X^{\vect,\gamma_j;\vecb,b,u,n+i+1}_{\tau_i}|\} \\
%&\quad+ |X^{\vect,\gamma_j;\vecb,b,u,n+i}_{\gamma}-X^{\vect,\gamma;\vecb,b,u,n+i}_{\gamma}|.
%\end{align*}
%Proceeding iteratively we find that
\begin{align*}
|X^{\vect,\gamma_j;\vecb,b,u}_\gamma-X^{\vect,\gamma;\vecb,b,u}_\gamma|
&\leq C(k+1)(\gamma-\gamma_j)+\sum_{i=1}^{k}\ett_{[\tau_i\leq \gamma]}\{|X^{\vect,\gamma_j;\vecb,b,u,n+i}_{\gamma} - X^{\vect,\gamma_j;\vecb,b,u,n+i}_{\tau_i}|
\\
&\quad+ |X^{\vect,\gamma_j;\vecb,b,u,n+i+1}_{\gamma} - X^{\vect,\gamma_j;\vecb,b,u,n+i+1}_{\tau_i}|\}+|X^{\vect;\vecb}_{\gamma}-X^{\vect;\vecb}_{\gamma_j}|
\\
&+|X^{\vect,\gamma_j;\vecb,b}_{\gamma}-X^{\vect,\gamma_j;\vecb,b}_{\gamma_j}|.
%
%+\sum_{i=0}^{k}\{|X^{\vect,\gamma_j;\vecb,b,u,n+i-1}_{\tau_i\wedge\gamma}-X^{\vect,\gamma_j;\vecb,b,u,n+i-1}_{\gamma_j}|
%\\
%&\quad+|X^{\vect,\gamma_j;\vecb,b,u,n+i}_\gamma-X^{\vect,\gamma_j;\vecb,b,u,n+i}_{\tau_i\wedge\gamma}|\},
\end{align*}
Furthermore, by H\"older's inequality we have
\begin{align*}
&\E\Big[\int_{0}^{\gamma}|X^{\vect,\gamma;\vecb,b,u}_{s}-X^{\vect,\gamma_j;\vecb,b,u}_{s}|^4ds]
\leq C\E[\gamma-\gamma_j]^{1/p}\E\big[((X^{\vect,\gamma;\vecb,b,u})^*_T+(X^{\vect,\gamma_j;\vecb,b,u})^*_T)^{4q}\big]^{1/q},
\end{align*}
where $\frac{1}{p}+\frac{1}{q}=1$. Now, by definition $\gamma$ is a predictable stopping time and the jump part of our SDDE is $\Prob$-a.s.~constant at predictable stopping times. We can, thus, apply Lemma~\ref{lem:SDDEcont} and the assertion follows.\qed

\begin{prop}
For all $(\vect,\vecb)\in\mcD^f$ and all $b\in \mcI^{-b_n}$, the process\\ $(\esssup_{u\in\mcU^k}V^{\vect,s\vee t_n;\vecb,b,u}_s:0\leq s\leq T)$ is in $\Sqlc^2$ for all $k\geq 0$.
\end{prop}

\noindent\emph{Proof.} Let $Y^{\vect;\vecb,k}_t:=\esssup_{u\in\mcU^k}V^{\vect;\vecb,u}_t$. To show that $Y^{\vect,\cdot\vee t_n;\vecb,b,k}_\cdot$ has a \cadlag version we consider
\begin{align*}
Y^{\vect,t';\vecb,b,k}_{t'}-Y^{\vect,t;\vecb,b,k}_{t}=(Y^{\vect,t';\vecb,b,k}_{t'}-Y^{\vect,t;\vecb,b,k}_{t'})+(Y^{\vect,t;\vecb,b,k}_{t'}-Y^{\vect,t;\vecb,b,k}_{t})
\end{align*}
where the second term on the right hand side goes to zero $\Prob$-a.s.~as $t'\searrow t$ by uniform integrability and right continuity of the filtration. Concerning the first term we have
\begin{align}
|Y^{\vect,t';\vecb,b,k}_{t'}-Y^{\vect,t;\vecb,b,k}_{t'}|
&\leq \sup_{u\in\mcU^k}\E\bigg[\int_t^T|f(s,X^{\vect,t';\vecb,b,u}_s) -f(s,X^{\vect,t;\vecb,b,u}_s)|ds
+|g(X^{\vect,t';\vecb,b,u}_T) - g(X^{\vect,t;\vecb,b,u}_T)|\nonumber
\\
&\quad +\sum_{j=1}^N|c_{\beta_{j-1},\beta_j}(\tau_j\vee t')-c_{\beta_{j-1},\beta_j}(\tau_j\vee t)|\Big|\mcF_{t'}\bigg]\nonumber
\\
&\leq \sup_{u\in\mcU^k}\E\bigg[\int_t^{t'}|f(s,X^{\vect,t';\vecb,b}_s) -f(s,X^{\vect,t;\vecb,b,u}_s)|ds\Big|\mcF_{t'}\bigg]\nonumber
\\
&\quad+k\sup_{s\in[t,t']}\sum_{b,b'\in\bar\mcI^2}|c_{b,b'}(t')-c_{b,b'}(s)|\nonumber
\\
&\quad +C(K)\sup_{u\in\mcU^k}\E\bigg[\int_{t'}^T|X^{\vect,t';\vecb,b,u}_s -X^{\vect,t;\vecb,b,u}_s|+|X^{\vect,t';\vecb,b,u}_T - X^{\vect,t;\vecb,b,u}_T|\Big|\mcF_{t'}\bigg]\nonumber
\\
&\quad+C\sup_{u\in\mcU^k}\E\Big[\sup_{r\in[t_n,T]}\ett_{[(X^{\vect,r;\vecb,b,u})_T^*\geq K]} (1+|(X^{\vect,r;\vecb,b,u})_T^*|^q) \Big|\mcF_{t'}\Big],\label{ekv:Ykcadlag}
\end{align}
for each $K>0$, by the local Lipschitz property of $f$ and $g$. Concerning the last term Doob's maximal inequality gives, for fixed $u\in\mcU^k$,
\begin{align*}
&\E\Big[\sup_{t\in[0,T]}\E\Big[\sup_{r\in[t_n,T]}\ett_{[(X^{\vect,r;\vecb,b,u})_T^*\geq K]} |(X^{\vect,r;\vecb,b,u})_T^*|^q \Big|\mcF_{t}\Big]^2\Big]
\leq C\E\Big[\sup_{r\in[t_n,T]}\ett_{[(X^{\vect,r;\vecb,b,u})_T^*\geq K]} |(X^{\vect,r;\vecb,b,u})_T^*|^{2q} \Big],
\end{align*}
Applying H\"older's inequality to the right hand side and taking the supremum over $\mcU$, we get
\begin{align*}
&\sup_{u\in\mcU}\E\Big[\sup_{t\in[0,T]}\E\Big[\sup_{r\in[t_n,T]}\ett_{[(X^{\vect,r;\vecb,b,u})_T^*\geq K]} |(X^{\vect,r;\vecb,b,u})_T^*|^q \Big|\mcF_{t}\Big]^2\Big]
\\
&\leq \sup_{u\in\mcU}(\Prob[\sup_{r\in[t_n,T]}(X^{\vect,r;\vecb,b,u})^*_T\geq K])^{1/2} \sup_{u\in\mcU}\Big(\E\Big[\sup_{r\in[t_n,T]}|(X^{\vect,r;\vecb,b,u})_T^*|^{4q}\Big]\Big)^{1/2}.
\end{align*}
Now, by Chebyshev's inequality and Proposition~\ref{prop:SDDEbndB},\\ $\sup_{u\in\mcU}\Prob[\sup_{r\in[t_n,T]}(X^{\vect,r;\vecb,b,u})_T^*\geq K]$ can be made arbitrarily small by choosing $K$ large. By monotonicity, it follows that the last term in \eqref{ekv:Ykcadlag} tends to zero, $\Prob$-a.s.~as $K\to\infty$. We conclude that $Y^{\vect,t';\vecb,b,k}_{t'}$ tends to $Y^{\vect,t;\vecb,b,k}_{t}$, $\Prob$-a.s.~when $t'\searrow t$ by right continuity of the switching costs in combination with Lemma~\ref{lem:xCADLAG} and it follows that $Y^{\vect,\cdot\vee t_n;\vecb,b,k}_\cdot$ has a \cadlag version.\\

Arguing as above we have that
\begin{align*}
&Y^{\vect,\gamma_j\vee t_n;\vecb,b,k}_{\gamma_j}-Y^{\vect,\gamma\vee t_n;\vecb,b,k}_{\gamma}
=(Y^{\vect,\gamma_j\vee t_n;\vecb,b,k}_{\gamma_j} - Y^{\vect,\gamma\vee t_n;\vecb,b,k}_{\gamma_j})+(Y^{\vect,\gamma\vee t_n;\vecb,b,k}_{\gamma_j}-Y^{\vect,\gamma\vee t_n;\vecb,b,k}_{\gamma}).
\end{align*}
Letting $j\to\infty$ the last term tends to zero $\Prob$-a.s.~by uniform integrability and quasi-left continuity of the filtration. Concerning the first term we have (where we for notational convenience assume that $\gamma,\gamma_j\in\mcT_{t_n}$)
\begin{align*}
\E\big[|Y^{\vect,\gamma_j;\vecb,b,k}_{\gamma_j} - Y^{\vect,\gamma;\vecb,b,k}_{\gamma_j}|\big]
%\\
%&\leq \sup_{u\in\mcU^k}\E\bigg[\int_{\gamma_j}^T|f(s,X^{\vect,\gamma_j;\vecb,b,u}_s) -f(s,X^{\vect,\gamma;\vecb,b,u}_s)|ds+|g(X^{\vect,\gamma_j;\vecb,b,u}_T) - g(X^{\vect,\gamma;\vecb,b,u}_T)|
%\\
%&\quad +\sum_{j=1}^N|c_{\beta_{j-1},\beta_j}(\tau_j\vee \gamma_j)-c_{\beta_{j-1},\beta_j}(\tau_j\vee \gamma)|\Big|\mcF_{\gamma_j}\bigg]
&\leq \sup_{u\in\mcU^k}\E\bigg[\int_{\gamma_j}^{\gamma}|f(s,X^{\vect,\gamma_j;\vecb,b,u}_s) -f(s,X^{\vect,\gamma;\vecb,b}_s)|ds\bigg]
\\
&\quad+ k\sum_{b,b'\in\bar\mcI^2}\sup_{\tau\in\mcT_{\gamma_j}}\E\big[|c_{b,b'}(\tau)-c_{b,b'}(\tau\vee\gamma)|\big]
\\
&\quad +C(K)\sup_{u\in\mcU^k}\E\bigg[\int_{\gamma}^T|X^{\vect,\gamma_j;\vecb,b,u}_s -X^{\vect,\gamma;\vecb,b,u}_s|+|X^{\vect,\gamma_j;\vecb,b,u}_T - X^{\vect,\gamma;\vecb,b,u}_T|\bigg]
\\
&\quad+C\sup_{u\in\mcU^{k+1}}\E\Big[\ett_{[(X^{\vect;\vecb,u})_T^*\geq K]}(1+ |(X^{\vect;\vecb,u})_T^*|^q) \Big]
\end{align*}
where the right hand side can be made arbitrarily small by Lemma~\ref{lem:xLCE} and quasi-left continuity of the switching costs. We conclude that
\begin{align*}
\lim_{j\to\infty} \E\left[|Y^{\vect,\gamma_j\vee t_n;\vecb,b,k}_{\gamma_j}-Y^{\vect,\gamma\vee t_n;\vecb,b,k}_{\gamma}|\right]=0,
\end{align*}
which implies that $Y^{\vect,\gamma_j\vee t_n;\vecb,b,k}_{\gamma_j}\to Y^{\vect,\gamma\vee t_n;\vecb,b,k}_{\gamma}$ in probability. Now since $Y^{\vect,\cdot\vee t_n;\vecb,b,k}_{\cdot}$ has left limits it follows that $Y^{\vect,\gamma_j\vee t_n;\vecb,b,k}_{\gamma_j}\to Y^{\vect,\gamma\vee t_n;\vecb,b,k}_{\gamma}$, $\Prob$-a.s.~and we conclude that $Y^{\vect,\cdot\vee t_n;\vecb,b,k}_{\cdot}\in\Sqlc^2$.\qed

\bigskip

By the above results we conclude that an optimal control for the hydropower planning problem does exist (under the assumptions detailed in this section). With a few notable exceptions (see \eg~\cite{KobilaHT,KobilaSC} in the case of singular control problems and Chapter 7 in \cite{OksenSulemBok} for examples of solvable impulse control problems) finding explicit solutions to impulse control problems is difficult. Instead we often have to resort to numerical methods to approximate the optimal control. A plausible direction for obtaining numerical approximations of solutions to the hydropower operators problem would be to further develop the Monte Carlo technique originally proposed for optimal switching problems in~\cite{CarmLud} (and later analyzed in~\cite{RAid}) to obtain polynomial approximations of $Y^{\vect,\vecb}$. Another possibility would be to apply the Markov-Chain approximations for stochastic control problems of delay systems developed in \cite{KushnerMCdelay}. However, a thorough investigation of either direction is out of the scope of the present work and will be left as a topic of future research.

\bibliographystyle{plain}
\bibliography{ConDepDyn_ref}
\end{document}